# EXACT SOLUTION TO THE CHOW-ROBBINS GAME FOR ALMOST ALL n, USING THE CATALAN TRIANGLE

John H. Elton

June 3, 2023

**Abstract**

The payoff in the Chow-Robbins coin-tossing game is the proportion of heads when you stop. Knowing when to stop to maximize expectation was addressed by Chow and Robbins(1965), who proved there exist integers $k_n$ such that it is optimal to stop when head minus tails reaches this. Finding $k_n$ exactly was unsolved except for finitely many cases by computer. We show $k_n = \left\lceil \alpha\sqrt{n} - 1/2 + \frac{(-2\zeta(-1/2))\sqrt{\alpha}}{\sqrt{\pi}} n^{-1/4} \right\rceil$ for almost all n, where $\alpha$ is the Shepp-Walker constant. This comes from our estimate $\beta_n = \alpha\sqrt{n} - 1/2 + \frac{(-2\zeta(-1/2))\sqrt{\alpha}}{\sqrt{\pi}} n^{-1/4} + O\left(n^{-7/24}\right)$ of real numbers defined by Dvoretzky(1967) for a more general Value function which is continuous in its first argument and easier to analyze. An $O(n^{-1/4})$ dependence was conjectured by Christensen and Fischer(2022) based on numerical evidence. Our proof uses moments involving Catalan and Catalan triangle numbers which appear in a tree resulting from backward induction, and a generalized backward induction principle. It was motivated by an idea of Häggström and Wästlund(2013) to use backward induction of upper and lower Value bounds from a horizon, which they used numerically to settle a few cases. Christensen and Fischer, with much better bounds, settled many more cases. We use Skorohod's embedding to get simple upper and lower bounds from the Brownian analog; our upper bound is the one found by Christensen and Fischer in a different way. We use them first for many more examples, but the new idea is to use them algebraically in the tree, with feedback to get a sharper Value estimate near the border, to settle almost all n.

jelton@bellsouth.net; elton@math.gatech.edu

## 1. Introduction.

1.1. *Background.* The Chow-Robbins game (also known as the $S_n / n$ problem) is a classical optimal stopping problem that can be stated in the form of a simple coin-tossing game, for which the payoff is the proportion of heads when you stop. The goal is to maximize your expected payoff. It seems to have been first posed by Breiman [1] in 1964, in the context of testing E.S.P., but was first analyzed by Chow and Robbins [4] in 1965. Let $S_n = \sum_{i=1}^{n} X_i$, where the $X_i$ are independent $\pm 1$ valued mean zero random variables, representing heads or tails in tossing a fair coin; this is a symmetric random walk. The object is to find a stopping time $\tau$ which is optimal in the sense that $E\left[\frac{S_\tau}{\tau}\right] = \sup_T E\left[\frac{S_T}{T}\right]$, the sup taken over stopping times (assumed a.s. finite). Chow and Robbins proved the existence of integers $0 < k_1 \le k_2 \le \ldots$ such that the stopping time $\tau = \inf\{n : S_n \ge k_n\}$ is optimal; no formula was given for the $k_n$ .

Next, in 1967, Dvoretzky [6] found a representation of an optimal stopping time in terms of a more general payoff, or Value, function. Define the Value starting from initial "position" $(u,n)$ under stopping time $T$ as

$$V(u,n,T) = E\left[\frac{u+S_T}{n+T}\right],$$

where $u$ is a real number, $n$ is a non-negative integer, and $T$ is an integer-valued stopping time for the symmetric random walk; and define

$$V(u,n) = \sup_T V(u,n,T) .$$

In fact, Dvoretzky allowed the $X_i$ to be more generally i.i.d. of mean zero and finite variance, but our paper is only concerned with the coin-tossing case. We emphasize that $u$ is allowed to be real, not just integer, unlike what Chow and Robbins considered in their proofs. This turns out to be quite significant. Dvoretzky proved that $V$ is a continuous function of the first argument $u$ , and the equation $V(\beta_n, n) = \beta_n / n$ uniquely defines a strictly increasing sequence of positive real numbers $0 < \beta_1 < \beta_2 \ldots$ such that the stop rule $\tau(u,n) = \min\{ j : u + S_j \ge \beta_{n+j}\}$ is optimal in the sense that $V(u,n) = V(u,n,\tau(u,n))$ . In particular, $\tau(0,0) = \min\{ j : S_j \ge \beta_j\}$ is optimal for the Chow-Robbins game. For our development below, we work with the real numbers $\beta_n$ rather than the integers $k_n$ because we can approximate them by approximating the Value function in the equation they satisfy; we can use real analysis. Interestingly, our approximation of the real numbers $\beta_n$ will allow us to give an *exact* formula for $k_n$ for almost all $n$ , in the sense of natural density (also called asymptotic density), as we'll discuss later. Dvoretzky showed that $.32 < \beta_n / \sqrt{n} < 4.06$ for sufficiently large $n$ , and conjectured that $\beta_n / \sqrt{n}$ approaches a limit.

In 1969, Shepp[13], and independently Walker[14], found a simple exact optimal stop rule for the continuous-time Brownian motion analog, which allowed them to prove Dvoretzky's conjecture. Let $W(t)$ be standard Brownian motion (Wiener process), and following Shepp's notation, define

$$V_W(u,b,T) = E\left[\frac{u+W(T)}{b+T}\right],\ V_W(u,b) = \sup_T V_W(u,b,T),$$

where $u$ and $b$ are real numbers with $b>0$. $T$ is a stopping time means it is a non-negative random variable that does not anticipate the future, and the sup taken over stopping times for which the expectation exists. Let $\alpha$ be the unique real root of $\alpha = (1-\alpha^2)\int_0^\infty \exp\left(\lambda\alpha - \lambda^2/2\right)d\lambda$. Computation gives $\alpha = .83992...$. Let $\tau_\alpha = \min\left\{t : u+W(t) \ge \alpha\sqrt{b+t}\right\}$. They proved $\tau_\alpha$ is the a.s. unique optimal stopping time, so $V_W(u,b) = V_W(u,b,\tau_\alpha)$; and

(1.1) $V_W(u,b) = (1-\alpha^2)\int_0^\infty \exp\left(\lambda u - \lambda^2 b/2\right)d\lambda,\ u \le \alpha\sqrt{b}$ ; $V_W(u,b) = u/b, u > \alpha\sqrt{b}$.

In other words, starting at time $b$, it is optimal to stop when you hit the square root boundary $\alpha\sqrt{b+t}$. Using the invariance principle (see e.g.[ 1, pg. 281]), Shepp [13, pg. 1005-1006] used the Brownian motion result to show that the optimal stopping boundary for the random walk game is asymptotic to $\alpha\sqrt{n}$; that is, $\lim_{n\to\infty} \beta_n/\sqrt{n} = \alpha$. But that does not give a way of knowing if it is optimal to stop at any specific position $(u,n)$ in the Chow-Robbins game, when $u$ is an integer. Medina and Zeilberger[10] discuss this distinction, pointing out that, at the time of their article (2009), not even $k_8$ was known (they refer to it as $\beta_8$, but we are adopting the notation of the original papers). They give some numerical data about early positions, and some good insight into the difficulty encountered when trying to rigorously decide whether to stop or go, from a given position.

In 2013, Häggström and Wästlund [6] showed, with a clever idea and the help of computer calculations, how to finesse the difficulty discussed in [10], and actually decide in some "early" positions whether or not stopping is optimal. Let $d$ be the number of heads minus the number of tails after $n$ flips. They expressed their results in terms of the number of heads, but we will give equivalent statements using $d$, to align with the usual notation. They observed that for a given value of $d$, there is a value , say $n_s(d)$, such that if $n \le n_s(d)$, you should stop, otherwise you should keep playing. Stating the stop rule in terms of $n$ as a function of $d$ is advantageous for computation: you only have to store the square root as many things as doing it the other way. Using upper and lower bounds for the value at any position, and using backward induction from "way out" (a horizon), they computed, for $d$ between 1 and 25, numbers $n_1(d) < n_2(d)$ such that if $n \le n_1(d)$ you should stop, and if $n \ge n_2(d)$ you should go on. The idea is that as you work backward from the horizon, those numbers should pull closer together. For $d$ less than 12, and for several more $d$'s between 13 and 25, for their horizon they found $n_2(d) = n_1(d)+2$ (note that $d$ and $n$ have the same parity), in which case $n_1(d) = n_s(d)$ and they therefore established the rigorous optimal stop rule for that $d$. Shepp's asymptotic value for the stopping rule put in this formulation is $n_s(d) \cong d^2/\alpha^2$. This is not so accurate: If you look at just the few cases where Häggström and Wästlund actually find $n_s(d)$, you can already see that it is off about first order in $d$, and actually it appears that $n_s(d) \cong (d^2+d)/\alpha^2$. Solving for $d$ in terms of $n$, this is consistent with the suggestion by Lai, Yao and Aitsahlia[9, pg. 768], that the stopping boundary for $d$ in terms of $n$, should be $\beta_n = \alpha\sqrt{n} - 1/2 + o(1)$.

But this limited amount of numerical data is not able to suggest anything more. By having much better upper and lower bounds, it is possible to get MUCH more data. Christensen and Fischer [5] (2022) gave much better upper and lower bounds for the optimal stopping value $V$ for the random walk, and used it to numerically settle very many more cases. They found the stop rule for $n$ up to 489241, which corresponds to $d$ up to about 588. We used the same upper bound that they did (with a different proof), but with a different lower bound, and settled yet again very many more cases than in [5], to attempt to get more numerical insight, as described in section 1.3. This eventually led to using our bounds to prove the theoretical results which are the subject of this paper.

1.2. *Embedding the random walk in Brownian motion, and Value bounds.* Our proofs of the upper and lower bounds on $V$ use the classical embedding of the random walk in $W$ using first-exit times, due to A. V. Skorohod (see e.g. [2, pg. 293]), which make the results seem rather intuitive. The embedding idea is quite natural: simply sample the Brownian path each time it changes by $\pm 1$, and you get a version of the symmetric random walk. Formally the properties follow from the strong Markov property. Let $T_0 = 0, T_n = \min\{t > T_{n-1} : |W(t) - W(T_{n-1})| = 1\}, n = 1, 2, \ldots$ . Then $W(T_n), n = 0, 1, \ldots$ has the same distribution as the process $S_n, n = 0, 1, 2, \ldots$, and $(T_n - T_{n-1}, W(T_n) - W(T_{n-1}))$ , $n = 1, 2, \ldots$ is an i.i.d. sequence, and $E[T_n] = n$. It is also evident that $T_n - T_{n-1}$ and $W(T_n) - W(T_{n-1})$ are independent, since the exit boundaries $W(T_{n-1}) + 1, W(T_{n-1}) - 1$ are symmetric about $W(T_{n-1})$. So in fact, the collection of random variables $T_n - T_{n-1}, n = 1, 2, \ldots$ is independent from the collection $W(T_n) - W(T_{n-1}), n = 1, 2, \ldots$, and by telescoping, the sequence $T_1, T_2, \ldots, T_n, \ldots$ is independent from $W(T_1), W(T_2), \ldots, W(T_n), \ldots$ .

LEMMA 1.1 (Christensen-Fischer[5, Theorem 1 pg. 3 ]). $V(u,b) \le V_W(u,b)$.

PROOF. Their proof uses superharmonic functions, in a more general setting. We give a proof using the embedding idea, as a preliminary for using it in our lower bound proof. Let $n^*$ be a stopping time for $S_n$. This induces a stopping time $T^* = T_{n^*}$ on $W$. So

$$V_W(u,b) \ge E\left[\frac{u+W(T^*)}{b+T^*}\right] = \sum_{n=0}^{\infty} E\left[\frac{u+W(T_n)}{b+T_n}\middle| T^*=T_n\right] P(n^*=n)$$

$= \sum_{n=0}^{\infty} E\left[\frac{u+W(T_n)}{b+n}\frac{b+n}{b+T_n}\middle| T^*=T_n\right] P(n^*=n)$ . But $T_n$ is independent of $W(T_n)$, and is also independent of $I_{T^*=T_n}$ since the latter is a function of $W(T_1), W(T_2), ..., W(T_n)$ (since $T^*$ is a stopping time) and $T_n$ is independent of $W(T_1), W(T_2), ..., W(T_n), ...$ . Thus

$\sum_{n=0}^{\infty} E\left[\frac{u+W(T_n)}{b+n}\frac{b+n}{b+T_n}\middle| T^*=T_n\right] P(n^*=n) = \sum_{n=0}^{\infty} E\left[\frac{b+n}{b+T_n}\right] E\left[\frac{u+W(T_n)}{b+n}\middle| T^*=T_n\right] P(n^*=n)$. By Jensen's inequality, $E\left[\frac{b+n}{b+T_n}\right] \ge \frac{b+n}{E[b+T_n]} = 1$, so $V_W(u,b) \ge \sum_{n=0}^{\infty} E\left[\frac{u+W(T_n)}{b+n}\middle| T^*=T_n\right] P(n^*=n)$

$= \sum_{n=0}^{\infty} E\left[\frac{u+S_n}{b+n}\middle| n^*=n\right] P(n^*=n) = E\left[\frac{u+S_{n^*}}{b+n^*}\right] = V(u,b,n^*)$. This is true for all stopping times $n^*$ for the random walk, so $V_W(u,b) \ge V(u,b)$, and the result is proved. □

In retrospect, it is as one would think: the random walk is a just a sampling of the Brownian motion, so naturally it can't do any better. There is the little matter of different denominators in the payoff, but Jensen's inequality goes the right way for that.

For a lower bound, we have

LEMMA 1.2. $V(u,b) \ge V_W(u,b)\left(1 - \frac{5}{12b}\left(1 + \frac{1}{\sqrt{b}}\right)\right),\ b > 1600$.

We'll give a detailed proof of Lemma 1.2 in Section 4. But the idea for our proof of this Lemma is simple enough: For $u$ below the Brownian boundary, run the Brownian motion until it hits the nearest integer to $\alpha\sqrt{b+t} - u$. With that stop rule, $E\left[\frac{u+W(T)}{b+T}\right]$ is about the same as $V_W(u,b)$ because we are so close to the boundary; we'll quantify this using a Fundamental Wald Identity of Shepp to obtain $E\left[\frac{u+W(T)}{b+T}\right] \ge V_W(u,b)\left(1 - \frac{1}{4b}\left(1 + \frac{1}{\sqrt{b}}\right)\right)$. But $W(T)$ is an integer, so in terms of the embedded random walk process $S_n = W(T_n)$, $E\left[\frac{u+W(T)}{b+T}\right] = \sum_{n=0}^{\infty} E\left[\frac{u+S_n}{b+T_n}\middle| n^*=n\right] P(n^*=n)$.

Proceed as in the proof of Lemma 1.1. Jensen's inequality goes the wrong way this time; but knowing the moments of the random time differences of the embedding, we can show that

$E\left[\frac{b+n}{b+T_n}\right] \le \left(1+\frac{1}{6b}\right)$, so

$E\left[\frac{u+W(T)}{b+T}\right] \le \left(1+\frac{1}{6b}\right)\sum_{n=0}^{\infty} E\left[\frac{u+S_n}{b+n}\middle| n^* = n\right] P(n^* = n) = \left(1+\frac{1}{6b}\right)V(u,b)$, giving the Lemma.

How good are these bounds? The data shows, and we will prove later, that if integer $u$ is just a hair more than $\alpha\sqrt{b}-1/2$ , then $V(u,b) = u/b$, but $V_W(u,b) \cong \left(1+.25b^{-1}\right)u/b$, so the Brownian upper bound overshoots the true value by a relative error of $O(b^{-1})$ at some places, for arbitrarily large $b$. The same examples give $V(u,b) \cong V_W(u,b)\left(1-.25b^{-1}\right)$, so we can't expect a lower bound of the form $V_W(u,b)\left(1-cb^{-1}\right)$ to do better than this; Lemma 1.2 gets lower bound $V_W(u,b)\left(1-.42b^{-1}\right)$. Theorem 6.2, later in Section 6, gives a greatly improved approximation to $V$ when $u$ is near the boundary, showing that $V_W$ is off from the true $V$ by essentially a relative error of $.25b^{-1}$ at half-integer values below the boundary, and the true $V$ is approximately piecewise linear in between, when the distance below the boundary is not more than about $b^{1/12}$. Theorem 6.2 does not give specific values for the constants (though that could be done with enough pain), and it was not used in our numerical work. Christensen and Fischer also give a lower bound, but our simple formula is convenient for our numerical work, and more importantly, for the later proof of the theoretical main results.

1.3 *Numerical exploration and speculation*. Using these good upper and lower bounds, we carried out the Häggström-Wästlund method numerically from a horizon of $n = 10^9$, and found $n_2(d) = n_1(d) + 2$ for $d$ up to 7995, and almost all cases up to 20,000, so that $n_s(d) = n_1(d)$ is settled for those. Stated another way, it settles all cases where $n < 7995^2/\alpha^2 \cong 9\times 10^7$, and most cases where $n < 20000^2/\alpha^2 \cong 5.7\times 10^8$, over half a billion, going considerably beyond Christensen and Fischer. Our computations did not really take much computer time, and we could eventually go a lot further, but we got pleasantly sidetracked by discovering the theoretical argument that ends up deciding the correct rule for almost all $n$. Having the stop-go boundaries based on $d$ rather than $n$, following Häggström and Wästlund, makes the computer algorithm extremely efficient and suitable for dealing with very large numbers. The spreadsheet for the answer has only 20,000 rows, rather than a billion.

Figure 1 is a graph of $n_2(d) - \left(d^2 + d\right)/\alpha^2$ for $d$ up to 20,000, from the spreadsheet; for almost all of those cases, and for all cases up to $d = 7995$ , $n_2(d) = n_s(d)$ . It is quite compelling! It appears thick since it is oscillating with about amplitude 2 around a square root curve.

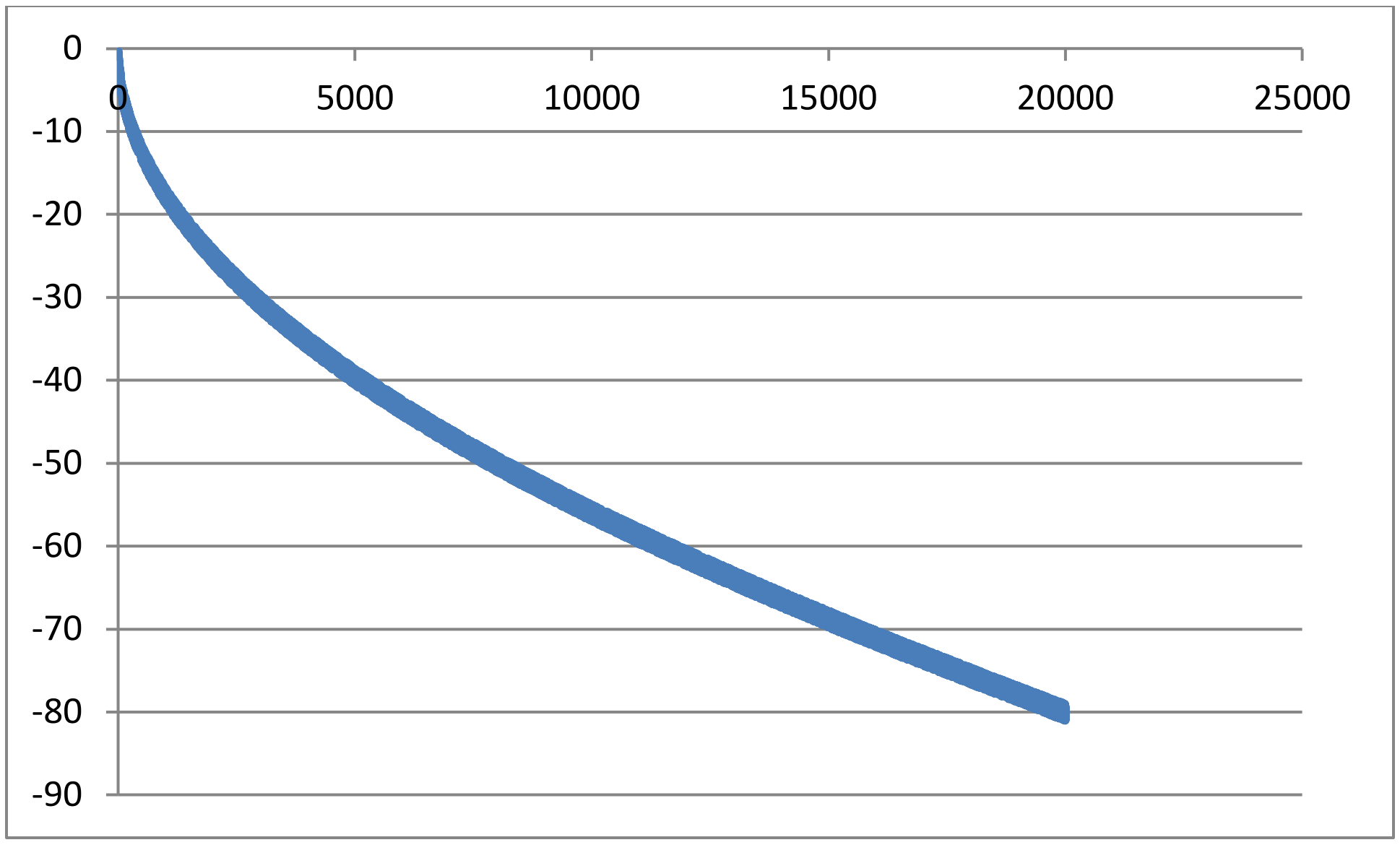

We originally speculated that $n_s(d) \cong \left(d^2+d\right)/\alpha^2 - \pi^{-1/2}\sqrt{d} \;\; + \varepsilon, \;\; -1 \le \varepsilon \le 1$ asymptotically, and it wiggles with oscillation amplitude about 1, probably according to number theory properties. But this coefficient of $\sqrt{d}$ turns out to be wrong, by about 1 percent. Theorem 1.3 will prove that the correct coefficient is $-\pi^{-1/2}\left(-4\zeta(-1/2)/\alpha\right)$. Since $-4\zeta(-1/2)/\alpha = .990...$ , nature had a laugh at us for jumping to conclusions!

Figure 2 is the detail for the 100 points at the large $d$ end of the curve; you can see the oscillation.

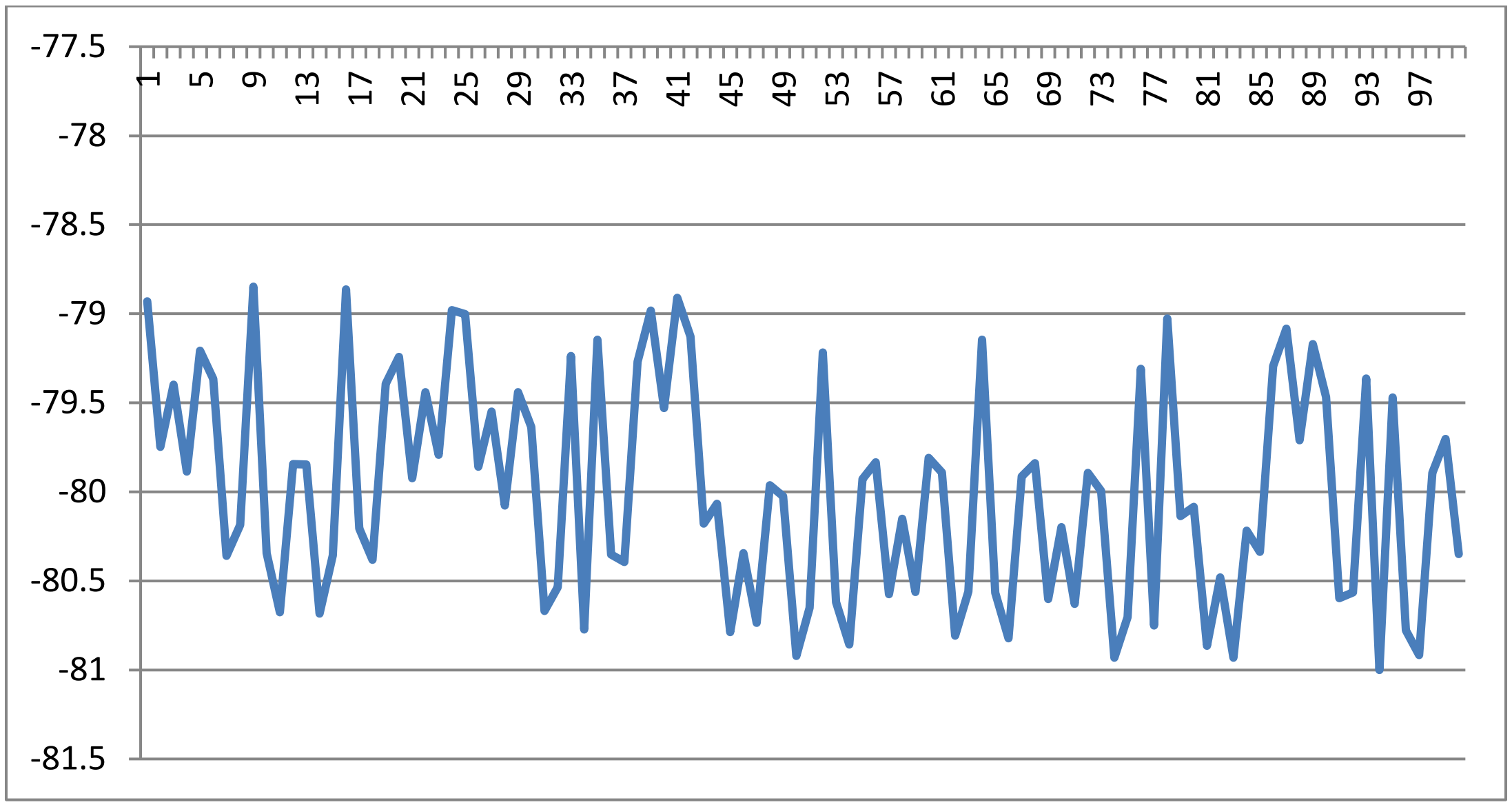

We now have, from theory, the correct coefficient for the $\sqrt{d}$ term, but there is still something suggested by the numerics that the theory has not yet reached. Assume $\varepsilon = O(1)$ . Using the binomial expansion, we can approximate the solution to $n=\left(d^2+d\right)/\alpha^2-c\sqrt{d}+\varepsilon$ for $d$, expressing the boundary in the more usual way with $d$ as a function of $n$. With a little algebra one gets

(1.2) Speculation. $\beta_n=\alpha\sqrt{n}-1/2+\dfrac{(-2\zeta(-1/2))\sqrt{\alpha}}{\sqrt{\pi}}n^{-1/4}+O(n^{-1/2})$? ?

Theorem 1.3 will prove that this coefficient of $n^{-1/4}$ is correct, but it will only get the error term to $O(n^{-7/24})$. There is still theoretical work to be done, to catch up with numerical speculation.

We remark that Christensen and Fischer had numerically observed an $n^{-1/4}$ dependence. They found that for $n$ up to 489241, the Chow-Robbins boundary is $k_n=\left\lceil \alpha\sqrt{n}-1/2+\dfrac{1}{7.9+4.54n^{1/4}}\right\rceil$ except for eight stray $n$'s in that range. Note that $\dfrac{1}{4.54}$ differs from $\dfrac{(-2\zeta(-1/2))\sqrt{\alpha}}{\sqrt{\pi}}$ by about 2.5 percent.

1.4. *Statement of main theorem*. Using the idea of Häggström and Wästlund to use backward induction from a horizon, but proceeding algebraically rather than numerically, we will be led to a tree with weights corresponding to Catalan numbers and Catalan triangle numbers, and a generalized backward induction principle. Before starting that development, we will up front state the main Theorem that eventually follows from it. As a Corollary, we give a formula for the optimal stopping rule for the Chow-Robbins game, for almost all $n$, in the sense of natural density. Note that previous results had established the exact stop rule for specific cases, up to some $n$, but now we are able to establish the *exact rule for almost all cases*, out to infinity. We remark that the appearance of the Riemann Zeta function $\zeta(-1/2)=-.207886...$ is not so mysterious; it appears because the analysis in Section 7 involves the asymptotic approximation of the sum of square roots of the first $k$ integers, the generalized harmonic number $H_k^{(-1/2)}$.

THEOREM 1. 3. $\beta_n=\alpha\sqrt{n}-1/2+\dfrac{(-2\zeta(-1/2))\sqrt{\alpha}}{\sqrt{\pi}}n^{-1/4}+O\left(n^{-7/24}\right)$.

THEOREM 1.4. $k_n=\left\lceil \alpha\sqrt{n}-1/2+\dfrac{(-2\zeta(-1/2))\sqrt{\alpha}}{\sqrt{\pi}}n^{-1/4}\right\rceil$ for almost all $n$.

Specifically, there exists $A>0$ and $n_0$ such that for $n\ge n_0$, the above equality holds if

$$\left\lceil \alpha\sqrt{n}-1/2+\frac{(-2\zeta(-1/2))\sqrt{\alpha}}{\sqrt{\pi}}n^{-1/4}-An^{-7/24}\right\rceil=\left\lceil \alpha\sqrt{n}-1/2+\frac{(-2\zeta(-1/2))\sqrt{\alpha}}{\sqrt{\pi}}n^{-1/4}+An^{-7/24}\right\rceil$$

PROOF. Figure 3 shows how Theorem 1.4 follows easily from Theorem 1.3. Let
$f_1(n)=\alpha\sqrt{n}\ -1/2+\dfrac{(-2\zeta(-1/2))\sqrt{\alpha}}{\sqrt{\pi}}n^{-1/4}-An^{-7/24}$ and
$f_2(n)=\alpha\sqrt{n}\ -1/2+\dfrac{(-2\zeta(-1/2))\sqrt{\alpha}}{\sqrt{\pi}}n^{-1/4}+An^{-7/24}$, with $A$ chosen according to Theorem 1.3 so that $f_1(n)<\beta_n<f_2(n)$ for all $n>n_0$. Let $u$ be a positive integer, and let $n_1,n_2,n_3$ satisfy $f_2(n_1)=u,\ f_1(n_2)=u,\ f_2(n_3)=u+1$, where for this purpose we extend the domain of $f_1,\ f_2$ so that $n_1,n_2,n_3$ are real, not necessarily integers. Assume also that $u$ is large enough so that $n_1>n_0$.

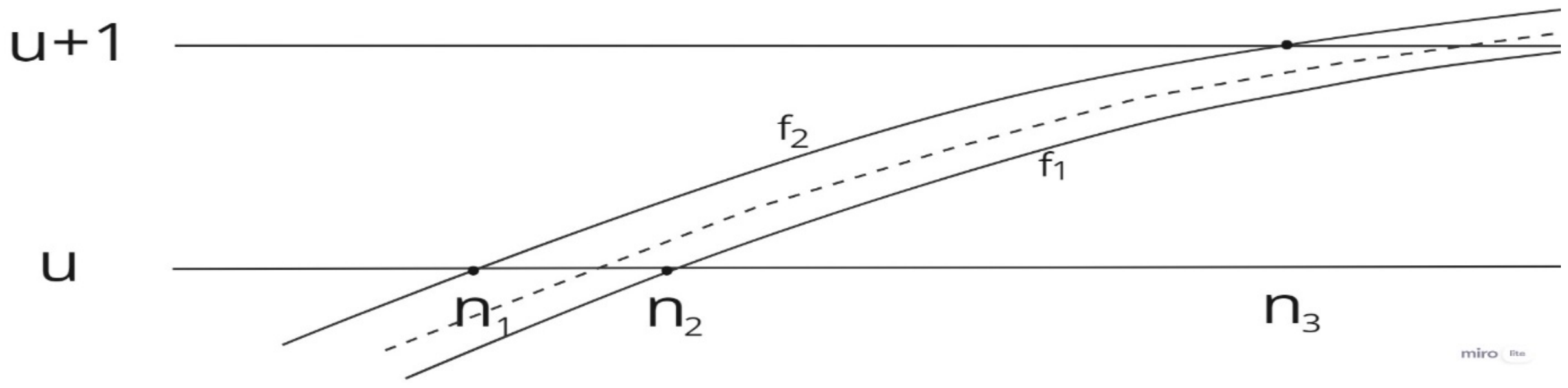

Graph of $f_1(n)$ and $f_2(n)$. $\beta_n$ is somewhere between them, shown dotted.

For $n_1\le n\le n_3$, $\lceil f_1(n)\rceil=\lceil f_2(n)\rceil$ implies $n_2\le n\le n_3$. For integers $n$ such that $n_2\le n\le n_3$, $u<\beta_n<u+1$, so it is optimal to stop at $u+1$, and for those integers
$k_n=u+1=\left\lceil \alpha\sqrt{n}\ -1/2+\dfrac{(-2\zeta(-1/2))\sqrt{\alpha}}{\sqrt{\pi}}n^{-1/4}\right\rceil$. Elementary estimates show $n_3-n_1\ge\sqrt{n_1}$, $n_2-n_1=O\left(n_1^{5/24}\right)$, so $(n_3-n_2)/(n_3-n_1)=1-O\left(n_1^{-7/24}\right)$. This proves the "almost all" assertion. □

The set where we are uncertain has a simple description as a union of intervals, such as $[n_1,n_2]$ in the picture above. The $i^{th}$ such interval is centered halfway between where $f_1$ and $f_2$ cross the horizontal line of height $i$, which is approximately at $i^2/\alpha^2$. The space between the $i^{th}$ and $(i+1)^{st}$ interval is $O(i)$. The length of the $i^{th}$ interval is $O\left(i^{5/12}\right)$. If our Speculation (1.2) were true, this length would be bounded.

We have not given a specific numerical value for $A$, though it could be done. We have resorted to expressing results in big-O notation, in spite of originally not wanting to do that, because we had wanted results that are usable for computer exploration. After frustration in some previous unrelated work exploring the Riemann Zeta (what else), where so many of the results in in the classic books are hopelessly unusable for computer exploration because of the long chains of big-O statements, we swore we would never do that. But we caved in and resorted to big-O, because of the onerous calculations in Sections 5, 6 and 7. O well (pun intended).

The goal of the rest of this paper is to prove Theorem 1.3, as well as to introduce the Catalan triangle for studying the problem. Having gotten this far, we are perhaps more optimistic than Medina and Zeilberger[10] about whether it is even possible to get a formula for $k_n$ for all $n$ (which might depend on the definition of "formula"). At least we are optimistic that Speculation (1.2) is conceivably provable with refinements of the techniques in this paper, or something similar. The proof of Theorem 1.3 will be given in Section 7, after quite a bit of preliminary development.

**2. Generalized Backward induction and the Catalan triangle. Plan for the proof of Theorem 1.3.** We start by patiently wading through some backward induction steps, rewarding us with a recognized pattern. We will want to decide whether to stop or continue when $u$ is below and near the Brownian boundary. Since we will be using the Brownian motion value function heavily for everything which follows, we switch back to using $(u,b)$ for position, as Shepp did. Shepp's wonderful paper [13] was our main source of inspiration and ideas early in our theoretical work on bounds.

We remark again that in our development from now on, $u$ is real, not just an integer, even though the Chow-Robbins game itself has only integer values for positions. This is important because our analysis depends on approximating $V(u,b)$, which Dvoretzky showed is continuous in $u$. The graphs in section 1.4 suggest the advantage of extending to the real case to do the analysis.

The famous *backward induction principle* of optimal stopping (see e.g. [3]) applied to this simple random walk is

(2.1) $$V(u,b)=\max\left\{\frac{u}{b},\frac{1}{2}V(u+1,b+1)+\frac{1}{2}V(u-1,b+1)\right\}.$$

Don't stop if $\frac{1}{2}V(u+1,b+1)+\frac{1}{2}V(u-1,b+1)>\frac{u}{b}$; stop otherwise.

This is the starting point for everything. Let's see how to use it, with a preliminary result. Numerical evidence showed that if $\delta=\alpha\sqrt{b}-u$ is larger than 1/2 minus a hair (the hair being on the order of $b^{-1/4}$), you should not stop. Our bounds on $V$ are not alone good enough to prove that without any backward steps, but let's see what we can get that way. Using only our lower bound from Lemma 1.2 (assuming $b>1600$), and a differential approximation $V_W(u,b)\ge\left(1+\delta^2b^{-1}\right)u/b$, from (3.6) of the next section, we get $V(u,b)\ge V_W(u,b)\left(1-.43b^{-1}\right)\ge\left(1+\delta^2b^{-1}\right)\left(1-.43b^{-1}\right)u/b$, and this is greater than $u/b$ if $\delta>.66$. Well, that's something: continue if $\delta>.66$. But just one step of backward induction with our lower bound will show how it begins to close in on 1/2. The distance of $u-1$ from the Brownian boundary is $\alpha\sqrt{b+1}-(u-1)=\alpha\sqrt{b}+h-u+1=\delta+1+h$, where $0<h<\alpha b^{-1/2}/2$.

So $\frac{1}{2}V(u+1,b+1)+\frac{1}{2}V(u-1,b+1)\ge\frac{1}{2}\frac{u+1}{b+1}+\frac{1}{2}V_W(u-1,b+1)\left(1-\frac{.43}{b+1}\right)$

$$\geq \frac{1}{2}\frac{u+1}{b+1}+\frac{1}{2}\frac{u-1}{b+1}\left(1+\frac{(1+\delta)^2}{b+1}\right)\left(1-\frac{.43}{b+1}\right)$$

$$=\frac{u}{b}-\frac{u}{b(b+1)}+\frac{1}{2}\frac{u\left(1-\frac{1}{\alpha\sqrt{b}-\delta}\right)}{(b+1)^2}\left((1+\delta)^2-.43-\frac{.43(1+\delta)^2}{b+1}\right).$$

For sufficiently large $b$, so that we can throw out small stuff, the condition for this to be greater than $u/b$ is clearly $(1+\delta)^2-.43>2$, or $\delta>.56$. But to be concrete, assume $b>1600$ and $u=\alpha\sqrt{b}-\delta>\alpha\sqrt{b}-.66$ (we already know to continue if $\delta$ greater than .66). Then some arithmetic shows that $\delta>.58$ is sufficient. That's an improvement. And we could go more steps back and get a better go bound.

Similarly, we can use our upper bound and close in on 1/2 from above, and it is convenient to do one step of that, to get a preliminary stop bound, to avoid annoyances later. Stop if $(V(u+1,b+1)+V(u-1,b+1))/2<u/b$. And $\alpha\sqrt{b+1}-(u-1)=\delta+1+h$ where $h$ is as before. Looking ahead to (3.5) of the next section, $V_W(u,b)\leq u/b+\alpha\delta^2 b^{-3/2}$. We have

$$\frac{1}{2}V(u+1,b+1)+\frac{1}{2}V(u-1,b+1)\leq\frac{1}{2}\frac{u+1}{b+1}+\frac{1}{2}\left(\frac{u-1}{b+1}+\alpha\frac{(1+\delta+h)^2}{(b+1)^{3/2}}\right)$$

$$=\frac{u}{b}-\frac{u}{b(b+1)}+\frac{1}{2}\alpha\frac{(1+\delta+h)^2}{(b+1)^{3/2}}=\frac{u}{b}-\frac{\alpha\sqrt{b}-\delta}{b(b+1)}+\frac{1}{2}\alpha\frac{(1+\delta+h)^2}{(b+1)^{3/2}}.$$

For sufficiently large $b$, the condition for this to be less than $u/b$ is clearly $(1+\delta)^2<2$, or $\delta<.414$. But again to be concrete, assuming $b>1600$, $\delta<.38$ can be shown to be sufficient.

With just one step of backward induction, we have now already narrowed the range for the boundary to:

LEMMA 2.1. For $b>1600$, $\alpha\sqrt{b}-.58<\beta_b<\alpha\sqrt{b}-.38$. Go if $\delta>.58$, stop if $\delta<.38$.

This will be useful later, so it is noted. But the goal is to get to $\alpha\sqrt{b}-1/2+cb^{-1/4}+o(b^{-1/4})$, by continuing way down the backward induction tree.

Without further ado, we proceed to systematize going backward. Continued backward induction leads to the following tree, where further branching is stopped at $u+1$, creating leaves of the tree at those nodes, which are boxed in the picture below. The $V$ value at a node of the tree is greater than or equal to the average of the $V$ values at its two parents. Figure 4 is a picture of the backward induction tree, showing 8 rows. The meaning of the coefficients (weights) will be explained shortly, though it is perhaps already obvious from the way backward induction works for this simple symmetric random walk.

**The Backward Induction Tree**

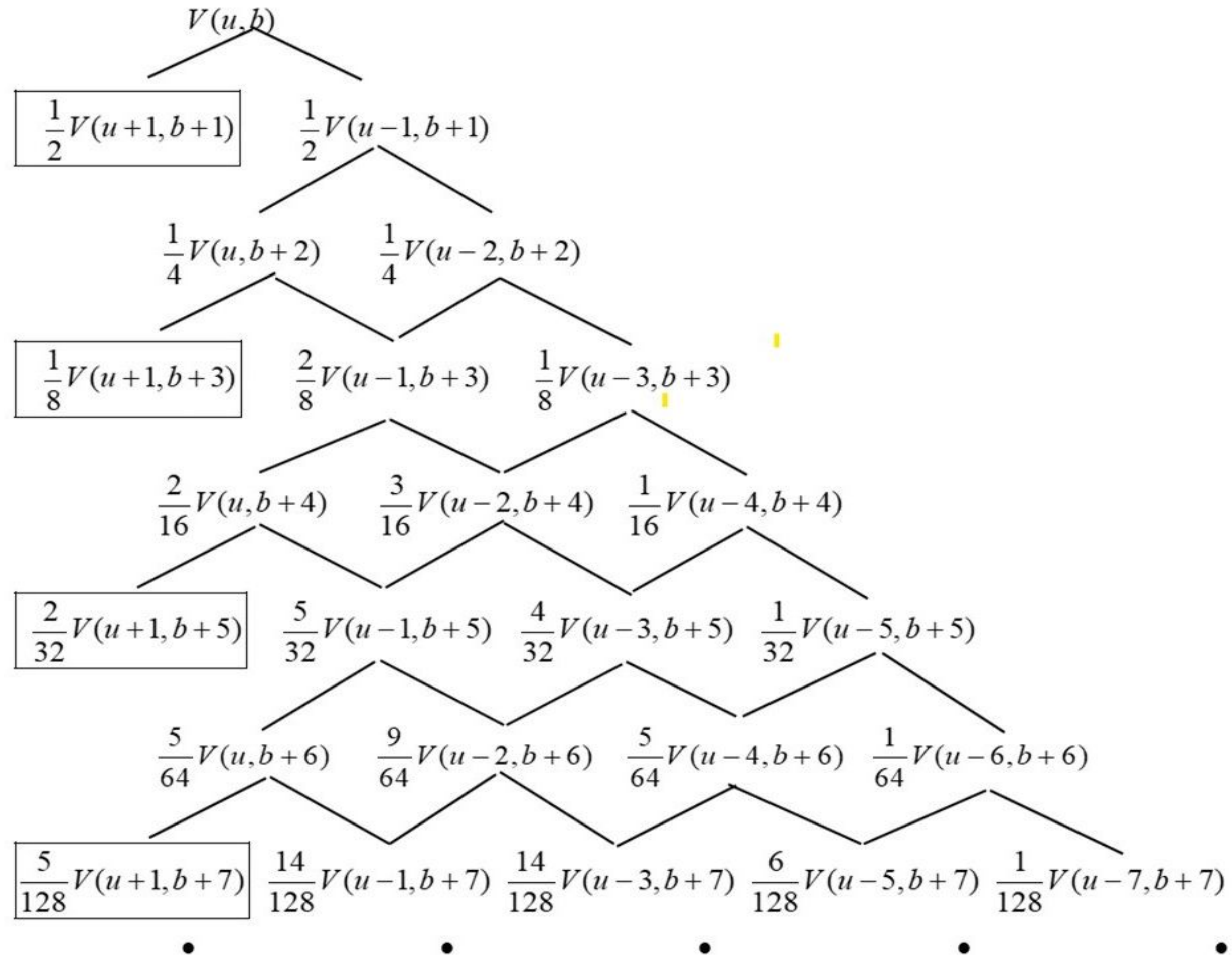

To explain the coefficients (the weights) displayed, consider this succession of inequalities, using only the basic backward induction inequality (2.1):

$V(u,b) \geq \frac{1}{2}V(u+1,b+1) + \frac{1}{2}V(u-1,b+1)$ ; $\frac{1}{2}V(u-1,b+1) \geq \frac{1}{4}V(u,b+2) + \frac{1}{4}V(u-2,b+2)$ ;

$$\frac{1}{4}V(u,b+2) + \frac{1}{4}V(u-2,b+2)$$
$$\geq \frac{1}{8}V(u+1,b+3) + \frac{1}{8}V(u-1,b+3) + \frac{1}{8}V(u-1,b+3) + \frac{1}{8}V(u-3,b+3)$$
$$= \frac{1}{8}V(u+1,b+3) + \frac{2}{8}V(u-1,b+3) + \frac{1}{8}V(u-3,b+3);$$

$$\frac{2}{8}V(u-1,b+3) + \frac{1}{8}V(u-3,b+3)$$
$$\geq \frac{2}{16}V(u,b+4) + \frac{2}{16}V(u-2,b+4) + \frac{1}{16}V(u-2,b+4) + \frac{1}{6}V(u-4,b+4)$$

$$= \frac{2}{16}V(u,b+4) + \frac{3}{16}V(u-2,b+4) + \frac{1}{6}V(u-4,b+4).$$

Etc. That's how the weights are generated, recursively. For any level of the tree, $V(u,b)$ is greater than or equal to the sum of the weights times values at the leaves at or above that level, plus the weights times values at the non-leaf nodes at that level. The picture shows the weights through level 7. Fortunately, we were patient enough to carry out the trivial calculations by hand through 7 rows, at which point the numbers were recognized: we recognize the leaf weight numerators as being the Catalan numbers, and the numerators of the weights at non-leaf nodes at a given level as being rows of a Catalan triangle, the Shapiro Catalan triangle (Shapiro [12]) (there are other things called "the" or "a" Catalan triangle in the literature, all related, but this is the relevant one for us).

Now formalize the notation and the recursion, which will prove that they are indeed those Catalan things. First, look at the nodes that are not leaves. Let $T(m, j) =$ coefficient of $2^{-m}V(u-j,b+m), m \geq 0, 0 \leq j \leq m$, and $T(m, j) = 0$ outside this range. The initial condition is $T(0,0) = 1$. The recursion that generates it is $T(m, j) = T(m-1, j-1) + T(m-1, j+1)$. Using the recursion produces the table in Fig. 5, through $m = 7$:

$\boldsymbol{T(m, j)}$

| m \ j | 0 | 1 | 2 | 3 | 4 | 5 | 6 | 7 |
|---|---|---|---|---|---|---|---|---|
| 0 | 1 | 0 | 0 | 0 | 0 | 0 | 0 | 0 |
| 1 | 0 | 1 | 0 | 0 | 0 | 0 | 0 | 0 |
| 2 | 1 | 0 | 1 | 0 | 0 | 0 | 0 | 0 |
| 3 | 0 | 2 | 0 | 1 | 0 | 0 | 0 | 0 |
| 4 | 2 | 0 | 3 | 0 | 1 | 0 | 0 | 0 |
| 5 | 0 | 5 | 0 | 4 | 0 | 1 | 0 | 0 |
| 6 | 5 | 0 | 9 | 0 | 5 | 0 | 1 | 0 |
| 7 | 0 | 14 | 0 | 14 | 0 | 6 | 0 | 1 |

In what follows we will only make use of the odd rows and columns; those correspond to the rows of the tree with leaves. Let $B(n,k) = T(2n-1, 2k-1), n \geq 1, k \geq 1$. The recursion for $T$ implies (in two steps) this recursion for $B$: $B(n,k) = B(n-1,k-1) + 2B(n-1,k) + B(n-1,k+1), n \geq 2, k \geq 2$, with initial conditions $B(1,1) = 1; B(1,k) = 0, k > 1$. This is the recursion and initial conditions for the Shapiro Catalan triangle, denoted in the literature by $B$; see Miana and Romero[11] for everything we will use about this Catalan triangle. Here are the first four rows:

$\boldsymbol{B(n,k)}$

| n \ k | 1 | 2 | 3 | 4 |
|---|---|---|---|---|
| 1 | 1 | 0 | 0 | 0 |
| 2 | 2 | 1 | 0 | 0 |
| 3 | 5 | 4 | 1 | 0 |
| 4 | 14 | 14 | 6 | 1 |

We now have a formula for the coefficients of the non-leaves in the odd rows of our backward induction tree in terms of well-known numbers:

(2.2) $\text{coefficient of } V(u-2j+1,b+2m-1) = 2^{-2m+1}B(m,j), m \ge 1, j \ge 1$.

There are simple formulas for the entries in $B$: $B(m,j) = \frac{j}{m}\binom{2m}{m-j} = \binom{2m-1}{m-j} - \binom{2m-1}{m-j-1}$.

Note $B(m,1) = \frac{1}{m}\binom{2m}{m-1} = \frac{1}{m+1}\binom{2m}{m} = \binom{2m-1}{m-1} - \binom{2m-1}{m-2} = C_m, m \ge 1$, the $m^{th}$ Catalan number. The zeroth Catalan number is defined by $C_0 = 1$.

For the leaves, let $L(n) =$ coefficient of $2^{-2n+1}V(u+1,b+2n-1), n \ge 1$. The first four are 1,1,2,5, and from the way the tree is built, it is seen that $L(n) = T(2n-2,0) = T(2n-3,1), n \ge 2$, with $L(1) = 1$. But $T(2n-3,1) = B(n-1,1) = C_{n-1}$ for $n \ge 2$. So, we have this formula for the leaf weights:

(2.3) $\text{coefficient of } V(u+1,b+2n-1) = 2^{-2n+1}C_{n-1}, n \ge 1$.

Returning to the inequality we started with: For any level, $V(u,b)$ is greater than or equal to the sum of the weights times values at the leaves at or above that level, plus the weights times values at the non-leaf nodes at that level. Looking at only the odd levels, in our general notation, this is, for $n \ge 1$,

$V(u,b) \ge \sum_{m=0}^{n-1} 2^{-2m-1}C_m V(u+1,b+2m+1) + 2^{-2n+1}\sum_{j=1}^{n} B(n,j)V(u-2j+1,b+2n-1)$. Define

(2.4) $TreeSum(n,u,b) = \sum_{m=0}^{n-1} 2^{-2m-1}C_m V(u+1,b+2m+1) + 2^{-2n+1}\sum_{j=1}^{n} B(n,j)V(u-2j+1,b+2n-1)$.

One property of $TreeSum(n,u,b)$ is immediate from the basic backward induction principle (2.1): $TreeSum(n+1,u,b) \le TreeSum(n,u,b),\ n \ge 1$, which implies $TreeSum(n,u,b) \le TreeSum(1,u,b)$ $= \left(V(u+1,b+1) + V(u-1,b+1)\right)/2$. Thus for any $n$, $TreeSum(n,u,b) > u/b \Rightarrow TreeSum(1,u,b) > u/b \Rightarrow V(u,b) > u/b$, so don't stop at $(u,b)$. The other direction is slightly more subtle. For $j \ge 0, k \ge 0$,

$V(u-j,b+k) = \sup_{n^*} E\left(\frac{u+S_{n^*}-j}{b+k+n^*}\right) = \sup_{n^*} E\left(\frac{u+S_{n^*}}{b+n^*}\left(1-\frac{k}{b+k+n^*}\right) - \frac{j}{b+k+n^*}\right)$

$\ge \sup_{n^*} E\left(\frac{u+S_{n^*}}{b+n^*}\left(1-\frac{k}{b+k}\right) - \frac{j}{b+k}\right) = \frac{V(u,b)b-j}{b+k}$, so $V(u,b) > \frac{u}{b} \Rightarrow V(u-j,b+k) > \frac{u-j}{b+k}$. In other words, if you shouldn't stop at $(u,b)$, then you shouldn't stop for a smaller $u$ or a larger $b$. Duh? This implies that if $TreeSum(1,u,b) > u/b$, then at every non-leaf node of the backward induction tree, the average of the values of its two children is greater than the ratio, so $TreeSum(n,u,b)$ $= TreeSum(1,u,b)$, and $TreeSum(n,u,b) > u/b$. This leads to

LEMMA 2.2. Extended backward induction principle. For any $n \geq 1$:
$V(u,b) = \max\{u/b, TreeSum(n,u,b)\}$; stop at $(u,b)$ if $TreeSum(n,u,b) \leq u/b$, else continue.
$TreeSum(1,u,b) > u/b \Rightarrow TreeSum(1,u,b) = TreeSum(n,u,b)$.

PROOF. From the previous paragraph, $TreeSum(n,u,b) > u/b \Rightarrow TreeSum(1,u,b) > u/b$ which implies $TreeSum(n,u,b) = TreeSum(1,u,b) = V(u,b)$. Now suppose $TreeSum(n,u,b) \leq u/b$. If $TreeSum(1,u,b) > u/b$, then by the previous paragraph, $TreeSum(n,u,b) = TreeSum(1,u,b) > u/b$, a contradiction. So $TreeSum(1,u,b) \leq u/b$, and $V(u,b) = u/b$ by definition, and you may stop. □

We can now explain the plan for the proof of Theorem 1.3, and give some indication of why it should work. The tree sum consists of a leaf sum, $S_L = \sum_{m=0}^{n-1} 2^{-2m-1} C_m V(u+1, b+2m+1)$ plus a row sum, $S_R = 2^{-2n+1} \sum_{j=1}^{n} B(n,j) V(u-2j+1, b+2n-1)$. The proof is accomplished in three stages.

Stage one is done in Section 5. For this stage, we consider a range of $n = c\sqrt{b}$ and $\delta = \alpha\sqrt{b} - u$ such that $V(u+1,b+2m+1) = (u+1)(b+2m+1)^{-1}$ exactly, and with a binomial expansion of that, the leaf sum can be approximated to any accuracy desired using simple formulas for $\sum_{m=0}^{n-1} 2^{-2m-1} m^k C_m$, which we use for $k = 0,1,2$. For the row sum, we use our simple upper and lower bounds for $V$ in terms of $V_W$, and approximate $V_W$ by a Taylor expansion about the boundary using four derivatives, and this leads to sums $2^{-2n+1} \sum_{j=1}^{n} j^k B(n,j)$, which have simple known formulas. What is the good of letting $n$ get big? The point is that the larger $n$, the more of the weight of the tree sum is on the leaves, for which the value is exactly known as long as $n$ does not get out of range, and the less is on the row sums, where we are limited in accuracy by our approximate bounds. This is in some way our algebraic manifestation of Häggström and Wästlund's idea to let the errors in the bounds wash out by moving the horizon back. In fact these sum formulas are all just simple expressions involving $n$ and the central binomial $2^{-2n}\binom{2n}{n} \cong (\pi n)^{-1/2}$, and $n$ will be of order $\sqrt{b}$, which hints at why $b^{-1/4}$ shows up in the answer. The first stage results in $O(b^{-1/4})$ upper and lower bounds on the stop rule; this is Lemma 5.3, which in fact was our original goal. But it is not able to get the exact coefficient of $b^{-1/4}$.

Stage two, in section 6, feeds the result of stage one back into the Value approximations developed in stage one, to obtain a much sharper approximation for $V$ near the boundary. It is perhaps the most conceptually tricky part of the proof, with a repeated feedback argument that will be better motivated when we get there. It shows that $V$ is essentially piecewise linear near the boundary, below and tangent at integer points of $\delta$, to the quadratic approximation to $V_W$ near the boundary. This is Theorem 6.2, perhaps of independent interest in showing the manner in which the Brownian Value overestimates $V$. Finally, stage three, in Section 7, uses this improved estimate of $V$ to go a bit further down the tree, by estimating the leaf values $V(u+1,b+2m+1)$ for a range of $m$ such that $V$ is no longer just the ratio. By going just far enough down the tree, we are able to get the upper and lower bounds to come together, within an $o(b^{-1/4})$ error, finding the exact coefficient of $b^{-1/4}$ and proving Theorem 1.3.

The plan is straightforward except perhaps for Section 6, and uses standard approximations, but lots of them, so it looks nasty when all the approximations and sums are written out; but it is not as bad as it looks. Like Mark Twain said about Wagner's music, it's not as bad as it sounds. (Love Wagner's music actually, but love the joke too, though Twain is not really its author).

**3. Approximating the Brownian motion value** $V_W$. From formula (1.1), one may differentiate under the integral sign to see that all the derivatives with respect to $u$ are positive for $u \le \alpha\sqrt{b}$, so they all take their maximum value on the boundary. To compute derivatives, and also for numerical work, it is best to write it in terms of standard functions. We have

$$V_W(u,b) = (1-\alpha^2)\int_0^\infty \exp\left(\lambda u - \frac{\lambda^2 b}{2}\right) d\lambda$$

$$= (1-\alpha^2) b^{-1/2} \exp\left(\frac{u^2}{2b}\right) \int_{-\infty}^{u/\sqrt{b}} \exp\left(-\frac{w^2}{2}\right) dw = (1-\alpha^2) b^{-1/2} G\left(u/\sqrt{b}\right)/g\left(u/\sqrt{b}\right)$$

$= (1-\alpha^2) b^{-1/2} H\left(u/\sqrt{b}\right)$, where $G$ and $g$ are the cdf and pdf of the standard normal distribution, respectively, and we have defined $H(x) = G(x)/g(x)$. This is related to the Mills Ratio by $H(x) = 1/g(x) - m(x)$, but we don't use that. On the boundary, $u_0 = \alpha\sqrt{b}$, $V_W(u_0,b) = u_0/b = \alpha b^{-1/2} = (1-\alpha^2) b^{-1/2} H\left(u_0/\sqrt{b}\right) = (1-\alpha^2) b^{-1/2} H(\alpha)$, so we get the equation that alpha satisfies as $\alpha = (1-\alpha^2) H(\alpha)$, which is useful for computing $\alpha$ using library function for the Normal, or converted to an Erf representation, to use that library function instead, if needed. We did that for the double-double precision numerical work for the extremely large number of positions we considered.

We will later use five derivatives of $V_W$ with respect to $u$ for the approximations. We'll use $D_u$ for derivative operator with respect to the first argument.

$D_u{}^n V_W(u,b) = (1-\alpha^2) b^{-(n+1)/2} H^{(n)}\left(u/\sqrt{b}\right)$, so we need to find the successive derivatives of $H$, which will be used also in Section 4. It satisfies $H'(x) = xH(x)+1$, which makes this straightforward, and can easily be made systematic in terms of polynomials in $x$, similar to Hermite polynomials. One shows

(3.1) $H^{(n)} = P_n H + Q_n$, where $P_{n+1} = P_n' + xP_n$, $Q_{n+1} = Q_n' + P_n$; $P_0 = 1$, $Q_0 = 0$.

So $D_u{}^n V_W(u,b) = (1-\alpha^2) b^{-(n+1)/2} \left(P_n\left(u/\sqrt{b}\right) H\left(u/\sqrt{b}\right) + Q_n\left(u/\sqrt{b}\right)\right)$.

REMARK. Just like for the Hermites, we have $P_n' = nP_{n-1}$ (an Appell sequence), so just like for the Hermites, there is a computationally practical recurrence, not involving derivatives: $P_{n+1}(x) = xP_n(x) + nP_{n-1}(x)$, the only difference from the Hermites being the positive sign. But we don't need to pursue this for our purposes, we just want a few derivatives.

On the boundary, $x = u_0/\sqrt{b} = \alpha\sqrt{b}/\sqrt{b} = \alpha$, and $H(\alpha) = \alpha(1-\alpha^2)^{-1}$, so

$D_u{}^n V_W(u_0,b) = (1-\alpha^2) b^{-n/2+1}\left(P_n(\alpha) H(\alpha) + Q_n(\alpha)\right) = b^{-(n+1)/2}\left(\alpha P_n(\alpha) + (1-\alpha^2) Q_n(\alpha)\right)$.

To get the polynomials, turn the recursion crank 5 times, resulting in

(3.2) $P_0 = 1,\ Q_0 = 0;\ P_1 = x,\ Q_1 = 1;\ P_2 = 1 + x^2, Q_2 = x;\ P_3 = 3x + x^3, Q_3 = 2 + x^2;$
$P_4 = 3 + 6x^2 + x^4, Q_4 = 5x + x^3;\ P_5 = 15x + 10x^3 + x^5, Q_5 = 8 + 9x^2 + x^4$.

In Section 4, we need the first 3 derivatives of $H$, which from (3.1) and (3.2) are

(3.3) $H^{(1)}(x) = xH(x) + 1, H^{(2)}(x) = (1+x^2)H(x) + x, H^{(3)}(x) = (3x + x^3)H(x) + 2 + x^2$.

The first 5 derivatives of $V_W$ on the boundary are

(3.4) $D_u V_W(u_0,b) = b^{-1}, D_u{}^2 V_W(u_0,b) = b^{-3/2} 2\alpha, D_u{}^3 V_W(u_0,b) = b^{-2}\left(2 + 2\alpha^2\right),$
$D_u{}^4 V_W(u_0,b) = b^{-5/2}\left(8\alpha + 2\alpha^3\right), D_u{}^5 V_W(u_0,b) = b^{-3}\left(8 + 16\alpha^2 + 2\alpha^4\right)$.

Let $\delta = u_0 - u = \alpha\sqrt{b} - u$. Using just the first two derivatives,

$V_W(u,b) = V_W(u_0,b) + (u-u_0) D_u V_W(u_0,b) + (u-u_0)^2 D_u{}^2 V_W(u_*,b)/2$

$\le \alpha b^{-1/2} - \delta b^{-1} + \delta^2 D_u{}^2 V_W(u_0,b)/2 = ub^{-1} + \alpha\delta^2 b^{-3/2}$. The inequality is because $D_u{}^2 V_W(u,b)$ is an increasing function (recall all derivatives are positive). In summary,

(3.5) $V_W(u,b) \le ub^{-1} + \alpha\delta^2 b^{-3/2}$.

Using the 3rd order term gives a lower bound that we used above in proving Lemma 2.1:

$V_W(u,b) \geq ub^{-1} + \alpha\delta^2 b^{-3/2} - (1+\alpha^2)\delta^3 b^{-2}/3 = ub^{-1} + \alpha b^{-1/2}\delta^2 b^{-2} - (1+\alpha^2)\delta^3 b^{-2}/3$

$= ub^{-1} + (u+\delta)\delta^2 b^{-2} - (1+\alpha^2)\delta^3 b^{-2}/3 = ub^{-1}\left(1+\delta^2 b^{-1}\right) + \delta^3 b^{-2}\left(1-(1+\alpha^2)/3\right)$

$\geq ub^{-1}\left(1+\delta^2 b^{-1}\right)$. Thus

(3.6) $V_W(u,b) \geq ub^{-1}\left(1+\delta^2 b^{-1}\right)$.

Using four derivatives gives the upper bound which will be used in proving Lemma 5.2:

(3.7) $V_W(u,b) \leq ub^{-1} + \alpha\delta^2 b^{-3/2} - (1+\alpha^2)\delta^3 b^{-2}/3 + (4\alpha+\alpha^3)\delta^4 b^{-5/2}/12$.

With the 5th derivative term we get the lower bound

(3.8) $V_W(u,b) \geq ub^{-1} + \alpha\delta^2 b^{-3/2} - (1+\alpha^2)\delta^3 b^{-2}/3 + (4\alpha+\alpha^3)\delta^4 b^{-5/2}/12 - (4+8\alpha^2+\alpha^4)\delta^5 b^{-3}/60$
.

**4. Proof of Lemma 1.2.** We'll need this estimate of the expected reciprocal which goes the other way.

LEMMA 4.1. $E\left[\frac{b+n}{b+T_n}\right] \leq \left[1+\frac{1}{6b}\right]$, where $T_n$ is as in section 1.2.

PROOF. We use these central moments: $E\left[\left(T_1-1\right)^2\right] = 2/3$, $E\left[(T_1-1)^3\right] = 16/15$; they can be found from the Laplace transform $f(t) = E[e^{-tT_1}] = \left(\cosh\sqrt{2t}\right)^{-1}, t>0$ [2, pg. 289]. Now

$$\frac{b+n}{b+T_n} = \frac{b+n}{b+n+T_n-n} = \left(1+\frac{T_n-n}{b+n}\right)^{-1} = 1 - \frac{T_n-n}{b+n} + \left(\frac{T_n-n}{b+n}\right)^2\left(1+\frac{T_n-n}{b+n}\right)^{-1}$$

$$= 1 - \frac{T_n-n}{b+n} + \frac{\left(T_n-n\right)^2}{\left(b+n\right)^2} - \frac{\left(T_n-n\right)^3}{\left(b+n\right)^2\left(b+n+T_n-n\right)}.$$

Note that $T_n - n = \sum_{j=1}^{n}(T_j - T_{j-1} - 1)$ is the sum of $n$ i.i.d. mean-zero random variables, so

$E\left[\left(T_n-n\right)^2\right] = E\left[\left(\sum_{j=1}^{n}(T_j-T_{j-1}-1)\right)^2\right] = \sum_{j=1}^{n} E\left[\left(T_j-T_{j-1}-1\right)^2\right] = nE\left[\left(T_1-1\right)^2\right] = 2n/3$, and

similarly, $E\left[\left(T_n-n\right)^3\right] = nE\left[\left(T_1-1\right)^3\right] = 16n/15$. In both cases, the cross terms disappear because of independence and mean-zero. Note the function $f(x)=x^3/(c+x)$ is convex for $x\ge -3c/2$, and $T_n-n\ge -3(b+n)/2$ is always true because $T_n\ge 0$. By Jensen,

$E\left[\frac{(T_n-n)^3}{(b+n+T_n-n)}\right]\ge \frac{E\left[(T_n-n)^3\right]}{(b+n+E[T_n-n])}=\frac{n}{b+n}>0$. Thus

$E\left[\frac{b+n}{b+T_n}\right]\le 1-\frac{E[T_n-n]}{b+n}+\frac{E\left[(T_n-n)^2\right]}{(b+n)^2}=1+\frac{(2/3)n}{(b+n)^2}$. It is easy to show that $n(b+n)^{-2}\le \frac{1}{4b}$ (worst case when $n=b$), so $E\left[\frac{b+n}{b+T_n}\right]\le 1+\frac{1}{6b}$. □

We need to estimate the loss in Value for the Brownian motion case if we get as close to the optimal boundary as possible while restricted to only stopping at integer steps from the start of the Brownian motion. Let $b$ be an integer, and $u<\alpha\sqrt{b}$. Let $f(t)=\left\lfloor \alpha\sqrt{b+t}\ -u+\ 1/2\right\rfloor$. Let $T$ be the first time $t$ that $W(t)=f(t)$ ( set $T=\infty$ if there is no such $t$). Let $F(t)=P[T\le t]$. We'll follow Shepp [13], using his Wald fundamental identity argument, except with this $f$.

LEMMA 4.2. For $\lambda\ge 0$, $\int_0^\infty dF(t)\exp\left\{\lambda\left(\left\lfloor \alpha\sqrt{b+t}\ -u+\ 1/2\right\rfloor\right)-\lambda^2 t/2\right\}=1$.

PROOF. This follows immediately from [13], Theorem 2, pg. 996, and the bottom of page 996. That theorem stipulated that $f$ be continuous, but the continuity was used only to justify the implication $T\ge t\Rightarrow W(t)\le f(t)$, and this is true when $f$ is monotone non-decreasing, without requiring continuity. □

PROOF OF LEMMA 1.2. Following Shepp, multiply both sides in Lemma 4.2 by $\alpha\exp\left\{\lambda u-\lambda^2 b/2\right\}$ and integrate over $\lambda$ from 0 to $\infty$, getting

$I=\alpha\int_0^\infty dF(t)\int_0^\infty \exp\left\{\lambda\left(\left\lfloor \alpha\sqrt{b+t}\ -u+\ 1/2\right\rfloor+u\right)-\lambda^2(b+t)/2\right\}d\lambda=\alpha\int_0^\infty \exp\left\{\lambda u-\lambda^2 b/2\right\}d\lambda$.

Define $r(t)=\left\lfloor \alpha\sqrt{b+t}\ -u+\ 1/2\right\rfloor-\alpha\sqrt{b+t}+u$, so $\left\lfloor \alpha\sqrt{b+t}\ -u+\ 1/2\right\rfloor=\alpha\sqrt{b+t}-u+r(t)$, where $-1/2\le r(t)<1/2$. Let $\varepsilon=\varepsilon(t)=r(t)(b+t)^{-1/2}$, so

$I=\alpha\int_0^\infty dF(t)\int_0^\infty \exp\left\{\lambda\sqrt{b+t}\left(\alpha+\varepsilon\right)-\lambda^2(b+t)/2\right\}d\lambda$. Making the change of variable $w=\lambda\sqrt{b+t}-(\alpha+\varepsilon)$ to complete the square in the integral yields

$I=\int_0^\infty dF(t)\alpha(b+t)^{-1/2}\exp\left\{\left(\alpha+\varepsilon\right)^2/2\right\}\int_{-(\alpha+\varepsilon)}^\infty \exp\left\{-w^2/2\right\}dw$

$= \int_0^\infty dF(t)\alpha(b+t)^{-1/2} G(\alpha+\varepsilon)/g(\alpha+\varepsilon) = \int_0^\infty dF(t)\alpha(b+t)^{-1/2} H(\alpha+\varepsilon)$, where $G$ and $g$ are the cdf and pdf of the standard normal, and $H = G/g$ was defined in Section 3. To summarize so far,

(4.1) $$I = \int_0^\infty dF(t)\alpha(b+t)^{-1/2} H(\alpha+\varepsilon).$$

This is a perturbation of what it would be for the Brownian optimal boundary, for $\varepsilon = 0$. To approximate the perturbation, the derivatives of $H$ given in Section 3 will be useful.
$H(\alpha+\varepsilon) = H(\alpha) + \varepsilon H'(\alpha) + \varepsilon^2 H''(\alpha)/2 + \varepsilon^3 H'''(z)/6$ for some $z$ between $\alpha$ and $\alpha+\varepsilon$ (note $\varepsilon$ can be negative). From (3.3), and recalling $H(\alpha) = \alpha/(1-\alpha^2)$, using the terms through the second derivative gives $H(\alpha) + \varepsilon H'(\alpha) + \varepsilon^2 H''(\alpha)/2 = H(\alpha) + \varepsilon(\alpha H(\alpha)+1) + \varepsilon^2\left(\left(1+\alpha^2\right)H(\alpha)+\alpha\right)/2$
$= H(\alpha)\left\{1+\varepsilon\left(\alpha+(1-\alpha^2)/\alpha\right)+\varepsilon^2\left(\left(1+\alpha^2\right)+(1-\alpha^2)\right)/2\right\} = H(\alpha)\left\{1+\varepsilon/\alpha+\varepsilon^2\right\}$. For the third derivative terms, we simply want to bound it. Section 3 noted all derivatives of $H$ are positive, so for $\varepsilon \le 0$, $\varepsilon^3 H'''(z) \le 0$. For $\varepsilon > 0$, $\varepsilon^3 H'''(z) \le \varepsilon^3 H'''(\alpha+\varepsilon)$, and also, since $H'$ is increasing,
$H(\alpha+\varepsilon) \le H(\alpha) + \varepsilon H'(\alpha+\varepsilon) = H(\alpha) + \varepsilon\left((\alpha+\varepsilon)H(\alpha+\varepsilon)+\alpha+\varepsilon\right)$. Solving, $H(\alpha+\varepsilon)$
$\le H(\alpha)\left(\alpha+\varepsilon(1-\alpha^2)\right)\left(1-\varepsilon(\alpha+\varepsilon)\right) \le 1.0152 H(\alpha)$ assuming $b > 1600$ so that
$\varepsilon = r(t)(b+t)^{-1/2} \le 1/80$. Then from (3.3), $H'''(\alpha+\varepsilon)$
$\le \left\{\left(3(\alpha+\varepsilon)+(\alpha+\varepsilon)^3\right)(1.0152) + \left(2+(\alpha+\varepsilon)^2\right)/H(\alpha)\right\} H(\alpha) \le 4.2 H(\alpha)$. Thus
$H(\alpha+\varepsilon) = H(\alpha) + \varepsilon H'(\alpha) + \varepsilon^2 H''(\alpha)/2 + \varepsilon^3 H'''(z)/6$ $H(\alpha+\varepsilon) \le H(\alpha)\left\{1+\dfrac{\varepsilon}{\alpha}+\varepsilon^2+.7|\varepsilon|^3\right\}$

$$\le \frac{\alpha}{1-\alpha^2}\left\{\frac{\alpha\sqrt{b+t}+r(t)}{\alpha\sqrt{b+t}} + \frac{1}{4(b+t)} + \frac{.7}{8(b+t)^{3/2}}\right\}$$

$$\le \frac{\alpha}{1-\alpha^2}\,\frac{\alpha\sqrt{b+t}+r(t)}{\alpha\sqrt{b+t}}\left\{1+\frac{\alpha\sqrt{b+t}}{4(b+t)\left(\alpha\sqrt{b+t}-1/2\right)}\left(1+\frac{.35}{(b+t)^{1/2}}\right)\right\}$$

$$\le \frac{\sqrt{b+t}}{1-\alpha^2}\left[\frac{W(t)+u}{b+t}\right]\left\{1+\frac{1}{4b}\left(1-\frac{1}{2\alpha\sqrt{b}}\right)^{-1}\left(1+\frac{.35}{\sqrt{b}}\right)\right\}$$

$$\le \frac{\sqrt{b+t}}{1-\alpha^2}\left[\frac{W(t)+u}{b+t}\right]\left\{1+\frac{1}{4b}\left(1+\frac{1}{b^{1/2}}\right)\right\}$$, using $b > 1600$ in the last step.

Referring back to (4.1), $\alpha\int_0^\infty \exp\left\{\lambda u - \lambda^2 b/2\right\} d\lambda = I = \int_0^\infty dF(t)\dfrac{\alpha}{\sqrt{b+t}} H(\alpha+\varepsilon)$

$$\le \int_0^\infty dF(t)\frac{\alpha}{1-\alpha^2}\left[\frac{W(t)+u}{b+t}\right]\left\{1+\frac{1}{4b}\left(1+\frac{1}{b^{1/2}}\right)\right\} = \frac{\alpha}{1-\alpha^2} E\left[\frac{W(T)+u}{b+T}\right]\left\{1+\frac{1}{4b}\left(1+\frac{1}{b^{1/2}}\right)\right\}$$, or

(4.2) $$E\left[\frac{W(T)+u}{b+T}\right]\left\{1+\frac{1}{4b}\left(1+\frac{1}{b^{1/2}}\right)\right\}\geq(1-\alpha^2)\int_0^\infty \exp\left\{\lambda u-\lambda^2 b/2\right\}d\lambda = V_W(u,b)$$

Now embed the random walk in the Brownian motion as in Section 1.2, with $S_n = W(T_n)$; whenever $W(t)$ is an integer, $W(t)=W(T_n)$ for some $n$. For stop rule $T$, $W(T)$ only takes integer values. Let $n*$ be the first $n$ such that $W(T)=W(T_n)$.

$$E\left[\frac{W(T)+u}{b+T}\right]=\sum_{n=0}^{\infty}E\left[\frac{u+S_n}{b+n}\frac{b+n}{b+T_n}\middle| n^*=n\right]P(n^*=n)=\sum_{n=0}^{\infty}E\left[\frac{b+n}{b+T_n}\right]E\left[\frac{u+S_n}{b+n}\middle| n^*=n\right]P(n^*=n)$$

$\leq\left(1+\frac{1}{6b}\right)V(u,b)$, using Lemma 4.1. Combining with (4.2),

$$V(u,b)\geq V_W(u,b)\left(1-\frac{1}{4b}\left(1+\frac{1}{b^{1/2}}\right)\right)\left(1-\frac{1}{6b}\right)\geq V_W(u,b)\left(1-\frac{5}{12b}\left(1+\frac{1}{b^{1/2}}\right)\right).$$

□

REMARK. Numerical data suggests that $V(u,b)\geq V_W(u,b)\left(1-\frac{1}{4b}\left(1+\frac{1}{b^{1/2}}\right)\right)$. Now the inequality $n\left(b+n\right)^{-2}\leq 1/(4b)$ used in proving Lemma 4.1 is a gross overestimate when $u$ is near the boundary $\alpha\sqrt{b}$, because in that case the time to crossing will probably be for $n$ small compared to $b$. On the other hand, when $u$ is far below $\alpha\sqrt{b}$, our replacing $\frac{1}{4(b+t)}$ by $\frac{1}{4b}$ in steps leading up to (4.2) is a significant overestimate, because $t$ will be typically of order $b$ in order to get to the boundary. But we will not pursue these things, because reducing the 5/12 to 1/4 in Lemma 1.2 would only slightly improve the lower estimate in Lemma 5.2, and the real goal is to prove Theorem 1.3.

**5. Computing Tree Sums.** This section begins the heavy computation. We've resorted to big-O notation rather than getting specific constants, for which we already apologized earlier. We will use big-O notation in inequalities, since we deal with approximate upper and lower bounds. The meaning will probably be clear from the context, but to be certain, before we get started, we'll define our notation. Let $f, g$ be real-valued functions of real variable $b$, with $g(b)>0$ for all sufficiently large $b$. We define $f(b)\leq O\left(g(b)\right)$ to mean there exists $b_0$ and $K>0$ such that $f(b)\leq Kg(b), b\geq b_0$. Usually only $f(b)=O\left(g(b)\right)$ is defined in texts, which requires $\left|f(b)\right|\leq Kg(b), b\geq b_0$. That would not be convenient for our purposes. We'll be writing expressions like, for example, $f(b)\leq U(b)+O\left(b^{-2}\right)$, to mean there exists $b_0$ and $K>0$ such that $f(b)\leq U(b)+Kb^{-2}$ for all $b\geq b_0$. That's the same as $f(b)-U(b)\leq O\left(b^{-2}\right)$ in our definition. $U$ is an approximate upper bound for $f$, but $f(b)$ could be arbitrarily far below $U(b)$. This seems very natural for dealing with inequalities rather than equalities, but does not seem to standard notation. Similarly, still assuming $g(b)>0$ for all sufficiently large $b$, define $f(b)\ge O\left(g(b)\right)$ to mean there exists $b_0$ and $K>0$ such that $f(b)\ge -Kg(b), b\ge b_0$. For example, $f(b)\ge U(b)+O\left(b^{-2}\right)$ means there exists $b_0$ and $K>0$ such that $f(b)\ge U(b)-Kb^{-2}$ for all $b\ge b_0$.

Here is another way to describe the situation, using set language. $O\left(g(b)\right)=\left\{h:\text{ there exists } K>0 \text{ and } b_0 \text{ such that } \left|h(b)\right|\le Kg(b) \text{ for all } b\ge b_0\right\}$. Note that $h\in O\left(g(b)\right)$ implies $\left|h\right|$, $-\left|h\right|$, and $-h$ are in $O\left(g(b)\right)$, and $O\left(g(b)\right)=-O\left(g(b)\right)$. Now $f(b)\le O\left(g(b)\right)$ means there exists $h\in O\left(g(b)\right)$ such that $f(b)\le h(b)$ for all $b$; $f(b)=O\left(g(b)\right)$ means there exists $h\in O\left(g(b)\right)$ such that $f(b)=h(b)$ for all $b$; and $f(b)\ge O\left(g(b)\right)$ means there exists $h\in O\left(g(b)\right)$ such that $f(b)\ge h(b)$ for all $b$. That unifies the settings: for some $h\in O\left(g(b)\right)$, the inequality or equality is true when $h(b)$ replaces $O\left(g(b)\right)$ in the expression. That concludes our discussion about big-O notation in this section.

Throughout this section, for integer $b$ and real number $u$, let $\delta=\alpha\sqrt{b}-u$. If $\delta<.38$ or $\delta>.58$, we already know whether to stop or keep going, by Lemma 2.1. But we will also want to estimate tree sums for a somewhat larger delta, to estimate Leaf sums in proving Theorem 6.2. We assert that $\delta\le b^p$, where $0\le p\le 1/10$. We also assert that $n=O(b^{1/2+p})$, where, as above, $n$ indexes how far we go down the tree, looking at only the odd rows (the ones with leaves). We prove Lemma 5.2 with this generality. In our first application of Lemma 5.2 , to get preliminary stop bounds, $p$ will be zero. But in Section 6 we'll use Lemma 5.2 with larger delta to get an improved estimate of Value, and then in Section 7 with larger $n$, to finally prove Theorem 1.3.

Let $G_n=2^{-2n}\binom{2n}{n}$, for which it is well-known as a central binomial coefficient that $G_n\cong\left(\pi n\right)^{-1/2}$; to be more precise, $\left(\pi(n+1/2)\right)^{-1/2}\le G_n\le\left(\pi n\right)^{-1/2}$. The formulas that will be needed for the sums and moments of Catalan numbers and Catalan triangle rows can be simply expressed in terms of $G_n$.

LEMMA 5.1. Catalan number and Catalan Triangle number sums.

(5.1) $\sum_{j=0}^{n-1}C_j 2^{-2j-1}=1-G_n;\sum_{j=0}^{n-1}C_j 2^{-2j-1}j=nG_n+G_n-1;\sum_{j=0}^{n-1}C_j 2^{-2j-1}j^2=\frac{1}{3}n^2G_n-\frac{4}{3}nG_n-G_n+1.$

(5.2) $2^{-2n+1}\sum_{j=1}^{n}B(n,j)=G_n;\ 2^{-2n+1}\sum_{j=1}^{n}jB(n,j)=\frac{1}{2};\ 2^{-2n+1}\sum_{j=1}^{n}j^2B(n,j)=nG_n;$

$$2^{-2n+1}\sum_{j=1}^{n}j^3B(n,j)=\frac{3n-1}{4};\ 2^{-2n+1}\sum_{j=1}^{n}j^4B(n,j)=n(2n-1)G_n;2^{-2n+1}\sum_{j=1}^{n}j^5B(n,j)=\frac{15n(n-1)+2}{8}.$$

PROOF. The first statement in (5.1) is proved quickly by induction. The other two are easily proved by telescoping; we omit the details. All of (5.2) can be found in Miana-Romero [9, pg. 5-6]. □

When proving Lemma 5.3, we'll be using the ratio for the value, on all the leaves, to compute the leaf sum. For getting an upper bound we'll be restricting $n$ to be small enough so that the ratio is in fact the correct value on all leaves, and for the lower bound we'll be restricting $n$ to be small enough so that the ratio (which is always a lower bound) will not be far off from the correct value. But in proving the final Theorem 1.1, later in Section 7, $n$ will be large enough so that the values at some of the leaves will be significantly larger than the ratio, and we have to take that into account. It will be convenient to split the leaf sum into two parts. For general $u$ and $b$, define $V_E(u,b) = V(u,b) - u/b$, the amount the value exceeds the ratio. Let

$S_{LR}(n,u,b) = \sum_{m=0}^{n-1} 2^{-2m-1} C_m \frac{u+1}{b+2m+1}$, the contribution from the ratio, and

$S_{LE}(n,u,b) = \sum_{m=0}^{n-1} 2^{-2m-1} C_m V_E(u+1, b+2m+1)$, from the excess over the ratio. So $S_L = S_{LR} + S_{LE}$.

For now, we'll only estimate $S_{LR}$, and leave $S_{LE}$ for Section 6 .

$$S_{LR} = \sum_{m=0}^{n-1} 2^{-2m-1} C_m \frac{u+1}{b+2m+1} = \frac{u+1}{b} \sum_{m=0}^{n-1} 2^{-2m-1} C_m \left(1 + \frac{2m+1}{b}\right)^{-1}$$

$= \frac{u+1}{b} \sum_{m=0}^{n-1} 2^{-2m-1} C_m \left(1 - \frac{2m+1}{b} + \gamma(m)\frac{(2m+1)^2}{b^2}\right)$, where $1 - \frac{2m+1}{b} \le \gamma(m) \le 1$. From (5.1),

$\sum_{m=0}^{n-1} 2^{-2m-1} C_m (2m+1)^2 = 4\left(\frac{n^2-4n-3}{3} G_n + 1\right) + 4(nG_n + G_n - 1) + 1 - G_n = \frac{4}{3}n^2 G_n - \frac{4}{3} nG_n - G_n + 1$

. We have this upper bound on $S_{LR}$, using $\gamma(m) \le 1$ and (5.1):

$$S_{LR} \le \frac{u+1}{b} \left\{ \sum_{m=0}^{n-1} 2^{-2m-1} C_m - \frac{1}{b} \sum_{m=0}^{n-1} 2^{-2m-1} C_m (2m+1) + \frac{1}{b^2} \sum_{m=0}^{n-1} 2^{-2m-1} C_m (2m+1)^2 \right\}$$

$$= \frac{u+1}{b}(1-G_n) - \frac{u+1}{b^2}\left(2nG_n - 1 + G_n - \frac{4}{3}\frac{n^2}{b} G_n + \frac{4}{3}\frac{n}{b} G_n - \frac{1}{b} + \frac{G_n}{b}\right)$$

$$= \frac{u+1}{b}(1-G_n) - \frac{u+1}{b^2}\left(2nG_n - 1 + G_n - \frac{4}{3}\frac{n^2}{b} G_n + O(b^{-1/2})\right).$$

We used $n = O(b^{1/2+p})$, $\frac{n}{b} G_n = O(b^{-3/4+p/2}) = O((b^{-1/2})$ for our range of $p$. For a lower bound, we use $\gamma(m) \ge 1 - 2b^{-1/2+p}$. Then

$$\frac{u+1}{b} \sum_{m=0}^{n-1} 2^{-2m-1} C_m \gamma(m) \frac{(2m+1)^2}{b^2} \ge \left(1 - 2b^{-1/2+p}\right) \frac{u+1}{b^2} \left(\frac{4}{3}\frac{n^2 G_n}{b} - \frac{4}{3}\frac{nG_n}{b} - \frac{G_n}{b} + \frac{1}{b}\right)$$

$$= \frac{u+1}{b}(1-G_n) - \frac{u+1}{b^2}\left( 2nG_n - 1 + G_n - \frac{4}{3}\frac{n^2}{b}G_n + O(b^{-1/2}) \right)$$, since

$b^{-1/2+p}\frac{n^2G_n}{b} = b^{-1/2+p}O\left(\frac{b^{3(1/2+p)/2}}{b}\right) = O(b^{-3/4+5p/2}) = O(b^{-1/2})$ . The lower bound is the same as the upper to that order. Since $u = O(n) = O(b^{1/2+p})$, we can write this as

(5.3) $$S_{LR} = \frac{u+1}{b}(1-G_n) - \frac{u+1}{b^2}\left( 2nG_n - 1 + G_n - \frac{4}{3}\frac{n^2}{b}G_n \right) + O(b^{-2+p}) .$$

That takes care of the Leaves. Now we have to add in the row sums and bound them, using the bounds on $V$ in terms of $V_W$, and approximating $V_W$ with terms from a Taylor expansion about the boundary. $V(u-2j+1, b+2n-1) \le V_W(u-2j+1, b+2n-1)$, and the displacement of $u-2j+1$ from the Brownian boundary is $u-2j+1-\alpha\sqrt{b+2n-1} = u-2j+1-\alpha\sqrt{b}\left(1+\frac{2n-1}{b}\right)^{1/2}$

$$= u-2j+1-\alpha\sqrt{b}\left(1+\frac{2n-1}{2b}\left(1-\tau(n)\frac{2n-1}{b}\right)\right) = -\left[2j-\left(1-\delta-\psi(n)\right)\right]$$, where

$\psi(n) = \alpha\frac{n-1/2}{\sqrt{b}}\left(1-\tau(n)\frac{2n-1}{b}\right)$, with $0 < \tau(n) < 1/4$. Now $\psi(n) \le \alpha n/\sqrt{b}$, and $\psi(n) = \alpha n/\sqrt{b} - O(b^{-1/2+2p})$. Thus from (3.7),

(5.4) $V_W(u-2j+1, b+2n-1)$

$$\le \frac{u-2j+1}{b+2n-1} + \alpha\frac{\left[2j-\left(1-\delta-\psi(n)\right)\right]^2}{\left(b+2n-1\right)^{3/2}} - \frac{\left(1+\alpha^2\right)}{3}\frac{\left[2j-\left(1-\delta-\psi(n)\right)\right]^3}{\left(b+2n-1\right)^2} + \frac{\alpha\left(4+\alpha^2\right)}{12}\frac{\left[2j-\left(1-\delta-\psi(n)\right)\right]^4}{\left(b+2n-1\right)^{5/2}}$$

By (3.8), we get a lower bound for $V_W(u-2j+1, b+2n-1)$ by subtracting $\frac{\left(4+8\alpha^2+\alpha^4\right)}{60}\frac{\left[2j-\left(1-\delta-\psi(n)\right)\right]^5}{\left(b+2n-1\right)^3}$ from this. These terms are to be weighted by $2^{-2n+1}B(n, j)$ and summed for $j$ from 1 through $n$. But $2^{-2n+1}\sum_{j=1}^{n} j^5 B(n, j) = O\left(n^2\right)$, so the result with this last term would be $O\left(n^2b^{-3}\right) = O\left(b^{-2+2p}\right)$, so the upper and lower bounds of the sum are the same to that order, and only involve the terms up through the 4th power of $j$. To get a lower bound for row sum of $V$, we will simply multiply these $V_W$ sums by $1-\frac{5/12}{b+2n-1}\left(1+\frac{1}{\sqrt{b+2n-1}}\right) = 1-\frac{5}{12b}+O\left(b^{-3/2+p}\right)$.

We'll carry both cases along at the same time. It's a bit long, so we get the terms one at a time and sum,

using Lemma 5.1 throughout. Some of the detail in obtaining the $O$ bounds from the asserted bounds on $n$ and $\delta$ is left to the reader. Let

$$s_1 = 2^{-2n+1}\sum_{j=1}^{n} B(n,j)\frac{u-2j+1}{b+2n-1} = \frac{u+1}{b+2n-1}2^{-2n+1}\sum_{j=1}^{n} B(n,j) - \frac{2}{b+2n-1}2^{-2n+1}\sum_{j=1}^{n} jB(n,j)$$

$$= \frac{(u+1)G_n - 1}{b+2n-1} = \frac{(u+1)G_n - 1}{b}\left(1 - \frac{2n-1}{b} + \frac{(2n-1)^2}{b^2} + O(n^3 b^{-3})\right)$$

$$= \frac{(u+1)G_n - 1}{b} - \frac{(u+1)G_n}{b^2}\left(2n - 1 - \frac{4n^2}{b}\right) + \frac{2n}{b^2} + O(b^{-2+p}).$$

For use in getting a lower bound, let $s_1' = s_1\left(1 - \frac{5}{12b} + O\left(b^{-3/2+p}\right)\right)$

$$= \frac{(u+1)G_n - 1}{b} - \frac{(u+1)G_n}{b^2}\left(2n - \frac{7}{12} - \frac{4n^2}{b}\right) + \frac{2n}{b^2} + O(b^{-2+p}).$$

Combining these with (5.3), cancellations simplify and clarify things.

$$S_{LR} + s_1 = \frac{u}{b} - \frac{u+1}{b^2}\left(4nG_n - 1 - \frac{16}{3}\frac{n^2 G_n}{b}\right) + \frac{2n}{b^2} + O(b^{-2+p});$$

$$S_{LR} + s_1' = \frac{u}{b} - \frac{u+1}{b^2}\left(4nG_n - 1 + \frac{5}{12}G_n - \frac{16}{3}\frac{n^2 G_n}{b}\right) + \frac{2n}{b^2} + O(b^{-2+p}).$$

The goal is to decide if the tree sum is greater than $u/b$ or not, so separating out $u/b$ like this reduces it to deciding if the combined other terms are negative or positive. Looking ahead, the other terms to be added will not have the $u+1$ factor, so we replace that now with an expression in terms of $n$ and $\delta$, and write the deviation from $u/b$ with factor $\alpha b^{-3/2}$. We want $\delta$ as the variable.

$$S_{LR} + s_1 - \frac{u}{b} = -\frac{\alpha\sqrt{b} + 1 - \delta}{b^2}\left(4nG_n - 1 - \frac{16}{3}\frac{n^2 G_n}{b}\right) + \frac{2n}{b^2} + O(b^{-2+p})$$

$$= -\frac{\alpha}{b^{3/2}}\left(1 + \frac{1-\delta}{\alpha\sqrt{b}}\right)\left(4nG_n - 1 - \frac{16}{3}\frac{n^2 G_n}{b}\right) + \frac{\alpha}{b^{3/2}}\frac{2n}{\alpha\sqrt{b}} + O(b^{-2+p}).$$ Simplifying,

(5.5) $$S_{LR} + s_1 - \frac{u}{b} = \frac{\alpha}{b^{3/2}}\left\{\begin{matrix} -4nG_n + 1 + \dfrac{2}{\alpha}\dfrac{n}{\sqrt{b}} - \dfrac{4}{\alpha}\dfrac{nG_n}{\sqrt{b}} + \dfrac{16}{3}\dfrac{n^2 G_n}{b} \\ +2\delta\left(\dfrac{2}{\alpha}\dfrac{nG_n}{\sqrt{b}}\right) \end{matrix}\right\} + O(b^{-2+p});$$

$$S_{LR} + s_1' - \frac{u}{b} = S_{LR} + s_1 - \frac{u}{b} + \frac{\alpha}{b^{3/2}}\left\{-\frac{5}{12}G_n\right\}.$$

Let $h = 1 - \delta - \psi(n)$, to simplify the writing for the next three terms.

$$s_2 = \frac{\alpha}{\left(b+2n-1\right)^{3/2}} 2^{-2n+1} \sum_{j=1}^{n} B(n,j)\left[2j-\left(1-\delta-\psi(n)\right)\right]^2$$

$$= \frac{\alpha}{b^{3/2}}\left(1-\frac{3}{2}\frac{(2n-1)}{b}+O(b^{-1+2p})\right)2^{-2n+1}\sum_{j=1}^{n} B(n,j)\left[4j^2-4jh+h^2\right]$$

$$= \frac{\alpha}{b^{3/2}}\left(1-\frac{3n}{b}+O(b^{-1+2p})\right)\left[4nG_n-2h+G_n h^2\right].$$

This can be arranged to (recall $\psi(n)=\alpha n/\sqrt{b}-O(b^{-1/2+2p})$)

(5.6 ) $$s_2 = \frac{\alpha}{b^{3/2}}\left\{\begin{array}{l} 4nG_n-2+2\alpha\frac{n}{\sqrt{b}}+G_n-\left(12-\alpha^2\right)\frac{n^2G_n}{b} \\ -2\alpha\frac{nG_n}{\sqrt{b}}+2\delta\left(1-G_n+\alpha\frac{nG_n}{\sqrt{b}}\right)+\delta^2G_n \end{array}\right\}+O(b^{-2+2p});$$

$$s_2' = s_2\left(1-\frac{5}{12b}+O\left(b^{-3/2+p}\right)\right)=s_2+O(b^{-2+2p}).$$

The correction factor for lower bound has no effect here.

For the next two terms, note $h=O\left(b^p\right)$ to help simplify.

$$s_3 = -\frac{(1+\alpha^2)}{3\left(b+2n-1\right)^2}2^{-2n+1}\sum_{j=1}^{n} B(n,j)\left[2j-\left(1-\delta-\psi(n)\right)\right]^3$$

$$= -\frac{(1+\alpha^2)}{3b^2}\left(1-\frac{4n}{b}+O(b^{-1+2p})\right)2^{-2n+1}\sum_{j=1}^{n} B(n,j)\left[8j^3-12j^2h+6jh^2-h^3\right]$$

$$= -\frac{\alpha}{b^{3/2}}\frac{(1+\alpha^2)}{3\alpha\sqrt{b}}\left(1+O(b^{-1/2+p})\right)\left[6n-2-12nG_nh+3h^2-G_nh^3\right].$$

The $h^2$ and $h^3$ terms contribute only $O(b^{-2+2p})$, so this becomes

(5.7) $$s_3 = \frac{\alpha}{b^{3/2}}\frac{(1+\alpha^2)}{\alpha}\left[-2\frac{n}{\sqrt{b}}+4\frac{nG_n}{\sqrt{b}}-4\alpha\frac{n^2G_n}{b}-2\delta\left(2\frac{nG_n}{\sqrt{b}}\right)\right]+O(b^{-2+2p});$$

$$s_3' = s_3\left(1-\frac{5}{12b}+O\left(b^{-3/2+p}\right)\right)=s_3+O(b^{-2+2p}).$$

$$s_4 = \frac{\alpha(4+\alpha^2)}{12\left(b+2n-1\right)^{5/2}}2^{-2n+1}\sum_{j=1}^{n} B(n,j)\left[2j-\left(1-\delta-\psi(n)\right)\right]^4$$

$$= \frac{\alpha(4+\alpha^2)}{12b^{5/2}}\left(1+O(b^{-1/2+p})\right)2^{-2n+1}\sum_{j=1}^{n} B(n,j)\left[\begin{array}{l}16j^4-32j^3h \\ +24j^2h^2-8jh^3+h^4\end{array}\right]$$

$$= \frac{\alpha}{b^{3/2}}\frac{(4+\alpha^2)}{12b}\left(1+O(b^{-1/2+p})\right)\left[\begin{array}{l}16n(2n-1)G_n-8(3n-1)h \\ +24nG_nh^2-4h^3+G_nh^4\end{array}\right].$$ This reduces to

(5.8) $s_4 = \frac{\alpha}{b^{3/2}}\left(\frac{8(4+\alpha^2)}{3}\frac{n^2 G_n}{b}\right) + O(b^{-2+2p})$;

$$s_4' = s_4\left(1 - \frac{5}{12b} + O\left(b^{-3/2+p}\right)\right) = s_4 + O(b^{-2+2p}).$$

Combining (5.5) through (5.8),

$$S_{LR} + s_1 + s_2 + s_3 + s_4 - \frac{u}{b}$$

$$= \frac{\alpha}{b^{3/2}}\left\{\begin{array}{l} -\left(1 - G_n - 2\alpha\frac{nG_n}{\sqrt{b}} + \frac{\alpha^2}{3}\frac{n^2 G_n}{b}\right) \\ +2\delta\left(1 - G_n - \alpha\frac{nG_n}{\sqrt{b}}\right) + \delta^2 G_n + O(b^{-1/2+2p}) \end{array}\right\};$$

$$S_{LR} + s_1' + s_2' + s_3' + s_4' - \frac{u}{b}$$

$$= \frac{\alpha}{b^{3/2}}\left\{\begin{array}{l} -\left(1 - \frac{7}{12}G_n - 2\alpha\frac{nG_n}{\sqrt{b}} + \frac{\alpha^2}{3}\frac{n^2 G_n}{b}\right) \\ +2\delta\left(1 - G_n - \alpha\frac{nG_n}{\sqrt{b}}\right) + \delta^2 G_n + O(b^{-1/2+2p}) \end{array}\right\}.$$

This proves

LEMMA 5.2. Bounds on Tree Sum. Let $0 \le p \le 1/10$, $n = O(b^{1/2+p})$ and $\delta \le b^p$, where $\delta = \alpha\sqrt{b} - u$. Then $\frac{\alpha}{b^{3/2}}\left\{C' + \delta B + \delta^2 A\right\} \le TreeSum(n,u,b) - S_{LE} - \frac{u}{b} \le \frac{\alpha}{b^{3/2}}\left\{C + \delta B + \delta^2 A\right\}$, where $C = -\left(1 - C_1 G_n\right) + O\left(b^{-1/2+2p}\right)$ and $C' = -\left(1 - C_1' G_n\right) + O\left(b^{-1/2+2p}\right)$, $B = 2\left(1 - B_1 G_n\right)$, $A = G_n$, with $C_1 = 1 + 2\frac{\alpha n}{\sqrt{b}} - \frac{1}{3}\frac{\alpha^2 n^2}{b}$, $C_1' = \frac{7}{12} + 2\frac{\alpha n}{\sqrt{b}} - \frac{1}{3}\frac{\alpha^2 n^2}{b}$, and $B_1 = 1 + \frac{\alpha n}{\sqrt{b}}$.

We can use this immediately to get preliminary $O(b^{-1/4})$ bounds on the stop rule. If we don't go too far down the tree, $S_{LE}$ will be zero. By Lemma 2.1, the value $V(u+1, b+2n-1)$ at the leaf will be just the ratio if $\alpha\sqrt{b+2n-1} - (u+1) \le .38$. But $\alpha\sqrt{b+2n-1} - (u+1)$ $= \alpha\sqrt{b}\left(1 + \frac{2n-1}{b}\right)^{1/2} - (\alpha\sqrt{b} - \delta + 1) \le \alpha\sqrt{b}\frac{1}{2}\frac{2n}{b} + \delta - 1 \le \frac{\alpha n}{\sqrt{b}} - .42$, assuming $\delta \le .58$ by Lemma 2.1. This will be less than .38 if $\frac{\alpha n}{\sqrt{b}} \le .8$, and then $S_{LE} = 0$. We make note of this for use later:

(5.9) $\frac{\alpha n}{\sqrt{b}} \le .8$ and $\delta \le .58$ implies $S_{LE} = 0$.

If $S_{LE} = 0$ and $C + \delta B + \delta^2 A \le 0$, the Tree Sum will not exceed the ratio, and if $C' + \delta B + \delta^2 A \ge 0$, the ratio will not exceed the Tree Sum. Let $\delta_0$ and $\delta_0^{'}$ be the positive numbers satisfying $C + \delta_0 B + \delta_0^{2} A = 0$ and $C' + \delta_0^{'} B + \delta_0^{'2} A = 0$, respectively. From Lemma 2.2, it follows that $\alpha\sqrt{b} - \delta_0^{'} \le \beta_b \le \alpha\sqrt{b} - \delta_0$ . In (a) of the next Lemma, we'll estimate $\delta_0$ and $\delta_0^{'}$ using some convenient $n's$ to get preliminary stop bounds from this. This was originally our goal and was a theorem, but we then discovered it was possible to prove the exact $b^{-1/4}$ term, and the following is now a Lemma. In addition, in (b) below we'll estimate $\delta_0$ and $\delta_0^{'}$ for a value of $n$ that will be used throughout Section 6 for convenience, even though it does not give quite as tight a bound as the one in (a).

LEMMA 5.3. Preliminary bounds on stop rule.

(a) $\alpha\sqrt{b} - 1/2 + \frac{.231}{\sqrt{\pi}} b^{-1/4} + O(b^{-1/2}) \le \beta_b \le \alpha\sqrt{b} - 1/2 + \frac{.429}{\sqrt{\pi}} b^{-1/4} + O(b^{-1/2})$.

(b) Let $n = \left\lfloor \frac{\sqrt{b}}{2} \right\rfloor$, and let $\delta_0$ and $\delta_0^{'}$ be the positive numbers satisfying $C + \delta_0 B + \delta_0^{2} A = 0$ and $C' + \delta_0^{'} B + \delta_0^{'2} A = 0$, respectively, where $A$, $B$, $C$, and $C'$ are as in Lemma 5.2. Then there exists $b_0$ such that $1/2 - \frac{.44}{\sqrt{\pi}} b^{-1/4} \le \delta_0 \le \delta_0^{'} \le 1/2 - \frac{.13}{\sqrt{\pi}} b^{-1/4},\ b \ge b_0$.

PROOF. We'll have $n = \Theta(b^{-1/2})$ in all cases here, so $A = G_n = O(n^{-1/2}) = O(b^{-1/4})$. For $\delta_0$, approximate the solution to the quadratic equation: $\delta_0 = \frac{-B + \sqrt{B^2 - 4AC}}{2A}$

$$= \frac{-B + B\left(1 - \frac{2AC}{B^2} - \frac{2A^2C^2}{B^4} - O(A^3)\right)}{2A} = -\frac{C}{B} - \frac{AC^2}{B^3} + O(b^{-1/2}),$$ since $A = O(b^{-1/4})$. Now

$$-\frac{C}{B} - \frac{AC^2}{B^3} = \frac{1 - C_1 G_n + O(b^{-1/2+2p})}{2(1 - B_1 G_n)} - G_n \frac{\left(1 - C_1 G_n + O(b^{-1/2+2p})\right)^2}{8\left(1 - B_1 G_n\right)^3}$$

$$= \frac{1}{2}\left(1 - C_1 G_n\right)\left(1 + B_1 G_n\right) - \frac{G_n}{8} + O(b^{-1/2+2p}) = \frac{1}{2} - \frac{1}{2}\left(C_1 - B_1 + \frac{1}{4}\right) G_n + O(b^{-1/2+2p}).$$ This gives

(5.10) $$\delta_0 = \frac{1}{2} - \frac{1}{2}\left(\frac{\alpha n}{\sqrt{b}} - \frac{1}{3}\frac{\alpha^2 n^2}{b} + \frac{1}{4}\right) G_n + O(b^{-1/2+2p}).$$

Changing $C$ to $C'$ in the above calculation gives

(5.11) $\delta_0' = \frac{1}{2} - \frac{1}{2}\left(\frac{\alpha n}{\sqrt{b}} - \frac{1}{3}\frac{\alpha^2 n^2}{b} - \frac{1}{6}\right)G_n + O(b^{-1/2+2p})$

To get an upper bound on the stop rule, we may assume $\delta \le .58$ by Lemma 2.1. Take $n = \left\lfloor .8\frac{\sqrt{b}}{\alpha} \right\rfloor$, so

$\frac{\alpha n}{\sqrt{b}} \le .8$ and $S_{LE} = 0$, and $\frac{\alpha n}{\sqrt{b}} = .8 + O(b^{-1/2})$, and $p = 0$. With this $n$,

$$A = G_n = \frac{1}{\sqrt{\pi n}}\left(1 + O(n^{-1})\right) = \sqrt{\frac{\alpha}{.8}}\frac{b^{-1/4}}{\sqrt{\pi}}\left(1 + O(b^{-1/2})\right).$$

From (5.10), $\delta_0 = \frac{1}{2} - \frac{1}{2}\left(.8 - \frac{1}{3}.64 + \frac{1}{4}\right)G_n + O(b^{-1/2}) = \frac{1}{2} - .4184\sqrt{\frac{\alpha}{.8}}\frac{b^{-1/4}}{\sqrt{\pi}} + O(b^{-1/2})$

$\ge \frac{1}{2} - \frac{.429}{\sqrt{\pi}}b^{-1/4} + O(b^{-1/2})$, which proves the upper bound.

For the lower bound, we decide to take $n = \left\lfloor 1.1\frac{\sqrt{b}}{\alpha} \right\rfloor$ so $\frac{\alpha n}{\sqrt{b}} = 1.1 + O(b^{-1/2})$, to get a little better result. For this $n$, it is not quite true that $S_{LE} = 0$, but the expression is still a lower bound. With that choice, (5.11) gives $\delta_0' = \frac{1}{2} - \frac{1}{2}\left(1.1 - \frac{1}{3}1.21 - \frac{1}{6}\right)\sqrt{\frac{\alpha}{1.1}}\frac{b^{-1/4}}{\sqrt{\pi}} + O(b^{-1/2}) \le \frac{1}{2} - \frac{.231}{\sqrt{\pi}}b^{-1/4} + O(b^{-1/2})$, completing the proof of (a).

For (b), we will take $n = \left\lfloor \frac{\sqrt{b}}{2} \right\rfloor$, and will make use of this in the next section, with larger $\delta$, where $p$ is not necessarily 0. Then (5.10) and (5.11) with this $n$ yield

$1/2 - \frac{.433}{\sqrt{\pi}}b^{-1/4} + O(b^{-1/2+2p}) \le \delta_0 \le \delta_0' \le 1/2 - \frac{.137}{\sqrt{\pi}}b^{-1/4} + O(b^{-1/2+2p})$, which for some $b_0$ implies the result in (b). □

## 6. Improving the estimate of V near the boundary, using feedback.

Our next goal is to use the above to get an improved estimate of $V$, which is very accurate when not too far from the boundary. This will allow us to go further down the backward induction tree by estimating the leaf values $V(u+1, b+2m-1)$ for larger $m$, for which $u+1$ is more than ½ below the boundary and the value is no longer just the ratio. This is used in Section 7 to obtain the correct coefficient for the $b^{-1/4}$ term. We shall show $V$ is approximately piecewise-linear near the boundary: that is Theorem 6.2 below. As a first step, Lemma 5.3 showed that $V(u,b) = \frac{u}{b}$ if $\delta \le 1/2 - \frac{.43}{\sqrt{\pi}}b^{-1/4}$ for $b$ sufficiently large, already an improvement over our previous upper bound $V(u,b)\le \frac{u}{b}+\frac{\alpha}{b^{3/2}}\delta^2$ coming from the Brownian motion case. The Brownian upper bound is too big by about $\frac{\alpha}{b^{3/2}}(1/4)$ when delta is near 1/2.

We introduce a bit more notation. In this section and the next, we'll be doing some estimates, in a chained fashion, an unbounded number of times, which could cause a problem if just using big-O notation. The book [7] gives a notation (attributed to de Bruijn) which will be convenient for us. Let $f(b)=L(g(b))$ mean that $|f(b)|\le|g(b)|$ for all $b$. In set language, $L(g)=\{f: |f(b)|\le|g(b)| \text{ for all } b\}$. Also, let $L^+(g)=\{f\in L(g): f(b)\ge 0 \text{ for all } b\}$. Similarly, let $O^+(g)$ be the non-negative members of $O(g)$.

We assert that throughout Section 6, $n=\left\lfloor \frac{\sqrt{b}}{2}\right\rfloor$. Lemma 5.3 (b) showed that using this $n$, $1/2-\frac{.44}{\sqrt{\pi}}b^{-1/4}\le\delta_0\le\delta_0^{'}\le 1/2-\frac{.13}{\sqrt{\pi}}b^{-1/4},\ b\ge b_0$, where $\delta_0$ satisfies $C+\delta_0B+\delta_0^{\,2}A=0$, and $\delta_0^{'}$ satisfies $C'+\delta_0^{'}B+\delta_0^{'2}A=0$. $A,B,C,C'$ are defined in Lemma 5.2. We also increase $b_0$ if necessary so that $b_0\ge 1000$, anticipating what is needed below. Now $n\ge\sqrt{b}/2-1=\sqrt{b}/2\left(1-2/\sqrt{b}\right)$, so $G_n\le\frac{1}{\sqrt{\pi n}}\le\sqrt{2.14/\pi}b^{-1/4}\le .83b^{-1/4}$ for $b\ge b_0$, a bound that will be used several times in the following.

In going further away from boundary, we'll use the following lemma. Recall that $V_E(u,b)=V(u,b)-u/b$.

LEMMA 6.1. For $1/2\le\delta\le b^p$, with $p\le 1/10$, and $n=\left\lfloor \frac{\sqrt{b}}{2}\right\rfloor$,

$$V_E(u,b)=\frac{\alpha}{b^{3/2}}\left(2(\delta-1/2)+L\left(\max\{\delta^2,1\}b^{-1/4}\right)\right)+S_{LE}(n,u,b), \text{ for } b\ge b_0.$$

PROOF. Let $1/2\le\delta\le b^p$. The Value is the same as the Tree Sum because $\delta\ge\delta_0^{'}$, so $u$ is below the stop boundary. $\delta_0^{'}=1/2-\gamma' b^{-1/4}$ for some $.13/\sqrt{\pi}\le\gamma'\le .44/\sqrt{\pi}$. Let $x=\delta-\delta_0^{'}=\delta-1/2+\gamma' b^{-1/4}$. From Lemma 5.2 in the lower bound case, $V_E(u,b)-S_{LE}$

$$\ge\frac{\alpha}{b^{3/2}}\left(C'+B\delta+A\delta^2\right)=\frac{\alpha}{b^{3/2}}\left(C'+B\delta_0^{'}+A\delta_0^{'2}+Bx+Ax^2+2Ax\delta_0^{'}\right)=\frac{\alpha}{b^{3/2}}\left(Bx+Ax^2+2Ax\delta_0^{'}\right).$$

Now $B=2(1-B_1G_n)=2-2\left(1+\alpha\frac{n}{\sqrt{b}}\right)G_n\ge 2-(2+\alpha)G_n$, and $A=G_n$, so $Bx+Ax^2+2Ax\delta_0^{'}$

$\geq 2(\delta - 1/2) + 2\gamma' b^{-1/4} + \left(-\left(2+\alpha\right)x + x^2 + 2x\delta_0'\right)G_n$

$= 2(\delta - 1/2) + 2\gamma' b^{-1/4} + x\left(x - \alpha - 1 + 2\gamma'\right)G_n$. But $x\left(x-\alpha-1+2\gamma'\right)$ has a minimum of $-\left(\alpha+1-2\gamma'\right)^2/4 \geq -.85$, so $x\left(x-\alpha-1+2\gamma'\right)G_n \geq -.85(.83)b^{-1/4} > -b^{-1/4}$, so

$V_E(u,b) - S_{LE} \geq \frac{\alpha}{b^{3/2}}\left(2(\delta-1/2) - b^{-1/4}\right)$.

Next, let $x = \delta - \delta_0 = \delta - 1/2 + \gamma b^{-1/4}$, where $.13/\sqrt{\pi} \leq \gamma \leq .44/\sqrt{\pi}$. From Lemma 5.2 in the upper bound case, similar to lower bound case we get $V_E(u,b) - S_{LE} \leq \frac{\alpha}{b^{3/2}}\left(Bx + Ax^2 + 2Ax\delta_0\right)$. Now $B \leq 2 - 2G_n$, so $Bx + Ax^2 + 2Ax\delta_0 \leq 2(\delta-1/2) + 2\gamma b^{-1/4} - 2xG_n + \left(x^2 + x\right)G_n$

$= 2(\delta-1/2) + 2\gamma b^{-1/4} + (x^2 - x)G_n$. If $\delta \leq 1$, then $x \leq 1$, so $2\gamma b^{-1/4} + (x^2-x)G_n \leq 2\gamma b^{-1/4} < b^{-1/4}$. If $\delta > 1$, $x^2 - x = \left(\delta - \left(1/2 - \gamma b^{-1/4}\right)\right)^2 - \delta + 1/2 - \gamma b^{-1/4} \leq \delta^2 - 1\left(1 - 2\gamma b^{-1/4}\right) + 1/4 - 1/2 - \gamma b^{-1/4}$

$\leq \delta^2 - 1$. So $2\gamma b^{-1/4} + (x^2-x)G_n \leq \left(.88/\sqrt{\pi}\right)b^{-1/4} + (\delta^2-1)\left(\sqrt{2.14}/\sqrt{\pi}\right)b^{-1/4} \leq \delta^2 b^{-1/4}$. □

We'll look at a few small ranges of delta first, to establish a pattern.
If $0 \leq \delta \leq 1/2$, then from (5.6), $S_{LE} = 0$. If $\delta \leq \delta_0$, $V(u,b) = \frac{u}{b}$. Now assume $\delta_0 \leq \delta \leq 1/2$. Let $x = \delta - \delta_0$, so $x \leq \gamma b^{-1/4}$. From Lemma 5.2, proceeding as in the above proof, $V_E(u,b) \leq \frac{\alpha}{b^{3/2}}\left(C + B\delta + A\delta^2\right) = \frac{\alpha}{b^{3/2}}\left(Bx + Ax^2 + 2Ax\delta_0\right)$. With this $x$, it is easy to see that $Bx + Ax^2 + 2Ax\delta_0 \leq b^{-1/4}$; we don't care to be any more precise than that. We get

(6.1) $0 \leq \delta \leq 1/2$ implies $V_E(u,b) = \frac{\alpha}{b^{3/2}}\left(L^+(b^{-1/4})\right)$, $b \geq b_0$ .

When $1/2 \leq \delta$, the value $V(u,b)$ is the tree sum, bigger than the ratio. The leaf sum part involves $V(u+1, b+2m+1)$, where $m \leq n-1$. The distance of $u+1$ below the boundary is $d = \max\left\{\alpha\sqrt{b+2m+1} - (u+1), 0\right\}$. But

$\alpha\sqrt{b+2m+1} - (u+1) = \alpha\sqrt{b}\left(1 + \frac{2m+1}{b}\right)^{1/2} - (\alpha\sqrt{b} - \delta - 1)$. We estimate this. Now

$1+\frac{m+1/2}{b}-\frac{(2m+1)^2}{8b^2}\le\left(1+\frac{2m+1}{b}\right)^{1/2}\le 1+\frac{m+1/2}{b}$, and $\frac{(2m+1)^2}{8b^2}\le\frac{\left(\sqrt{b}\right)^2}{8b^2}=\frac{1}{8b}$, so

$1+\frac{m}{b}\le\left(1+\frac{2m+1}{b}\right)^{1/2}\le 1+\frac{m}{b}+\frac{1}{2b}$, so $\frac{\alpha m}{\sqrt{b}}\le\alpha\sqrt{b}\left(1+\frac{2m+1}{b}\right)^{1/2}-\alpha\sqrt{b}\le\frac{\alpha m}{\sqrt{b}}+\frac{\alpha}{2\sqrt{b}}$. Thus

(6.2) $d=\left\{\alpha\sqrt{b+2m+1}-(u+1)\right\}^{+}=\left\{\frac{\alpha m}{\sqrt{b}}+\delta-1+L^{+}\left(\alpha b^{-1/2}/2\right)\right\}^{+}$, $m\le n-1$. We'll use this several places below.

Now consider $1/2\le\delta\le 1$. Since $\delta\le 1$, $d\le\frac{\alpha m}{\sqrt{b}}+L^{+}\left(\alpha b^{-1/2}/2\right)$. But $m\le\frac{\sqrt{b}}{2}$, so $d\le\frac{\alpha}{2}+L^{+}\left(\alpha b^{-1/2}/2\right)$, and computation shows this is less than $\frac{1}{2}-\frac{.44}{\sqrt{\pi}}b^{-1/4}$ if $b\ge b_0\ge 400$, assumed. So $V_E(u+1,b+2m+1)=0$ (that is, the value is just the ratio), and $S_{LE}=0$. Then Lemma 6.1 gives

(6.3) $1/2\le\delta\le 1$ implies $V_E(u,b)=\frac{\alpha}{b^{3/2}}\left(2(\delta-1/2)+L(b^{-1/4})\right)$, $b\ge b_0$.

Our improved estimate so far for $V(u,b)$, $0\le\delta\le 1$, is piecewise linear plus an $O(b^{-1/4})$ correction; zero for delta from 0 to 1/2, then slope 2 from 1/2 to 1. The piecewise linear part is accurate to order $b^{-1/4}$ (relative to $b^{-3/2}$). Compare this to our previous upper bound estimate, the Brownian value, $V_W(u,b)\le\frac{u}{b}+\frac{\alpha}{b^{3/2}}\delta^2$ and about equal to that for delta small compared to $\sqrt{b}$; that is, quadratic for delta not too big. Our piecewise linear function matches the quadratic one at $\delta=0$ and $\delta=1$, and is tangent to the parabola at those points. At $\delta=1/2$, the piecewise linear function is below the parabola by 1/4. That's where the Brownian upper bound is worst.

Now we show that same linear piece continues up to 3/2 . Let $1\le\delta\le 3/2$.
The distance of $u+1$ from the boundary $\alpha\sqrt{b+2m+1}$ is, from (6.2),
$d=\frac{\alpha m}{\sqrt{b}}+\delta-1+L^{+}\left(\alpha b^{-1/2}/2\right)\le\frac{\alpha}{2}+1/2+L^{+}\left(\alpha b^{-1/2}/2\right)\le 1$ for $b\ge b_0$ . If also $d\ge 1/2$, we can feed this into (6.3) with $d$ in the role of $\delta$ and $b+2m+1$ in the role of $b$. *This is the key feedback idea that* *will be used in proving* Theorem 6.1. If $d \le 1/2$, (6.3) still holds if we replace $d-1/2$ by $\{d-1/2\}^{+}$.

So $V_E(u+1,b+2m+1) = \frac{\alpha}{(b+2m+1)^{3/2}}\left(2\{d-1/2\}^{+} + L(b^{-1/4})\right)$.

Note $d-1/2 \le \frac{\alpha m}{\sqrt{b}} + \alpha b^{-1/2}/2$. Use (5.1) to sum:

$$S_{LE} \le \sum_{m=0}^{n-1} 2^{-2m-1} C_m \frac{\alpha}{(b+2m+1)^{3/2}}\left(2\frac{\alpha m}{\sqrt{b}} + \alpha b^{-1/2} + b^{-1/4}\right)$$

$\le \frac{\alpha}{b^{3/2}}\left(nG_n 2\frac{\alpha}{\sqrt{b}} + \alpha b^{-1/2} + b^{-1/4}\right) \le \frac{\alpha}{b^{3/2}}\left(2b^{-1/4}\right)$ for $b \ge b_0 \ge 1000$. Using Lemma 6.1, adding in this estimate of $S_{LE}$, gives the same answer as (6.3), except for increasing the error bound, which then covers (6.3) as well:

(6.4): $1/2 \le \delta \le 3/2$ implies $V_E(u,b) = \frac{\alpha}{b^{3/2}}\left(2(\delta-1/2) + L(5b^{-1/4})\right)$, $b \ge b_0$.

Now that we see that the piecewise linear function below and tangent to the quadratic at integers gets us this far, we guess that this pattern continues (up to some error), and Theorem 6.2 will prove this. By stepping along one 1/2 unit at a time, and feeding the result back into the previous step the way we did to get to 3/2, we will get our piecewise linear estimate that will be good enough for obtaining the right coefficient for the $b^{-1/4}$ term in the stop boundary.

THEOREM 6.2. Value near the boundary. Assume $b \ge b_0$. Then

$V_E(u,b) = \frac{\alpha}{b^{3/2}}\left(2j(\delta-1/2) - j(j-1) + L\left(M_j b^{-1/4}\right)\right)$, $j-1/2 \le \delta \le j+1/2, 1 \le j \le b^{1/10}$, where $j$ is integral, and $M_j \le 5j^3$. $M_1 = 5$, and $M_j$ satisfies the recursion $M_{j+1} = M_j + 4j^2 + 7j + 11, j \ge 1$.

Fig. 7 is a graph of the piecewise linear function $2j(\delta-1/2) - j(j-1)$ for $0 \le \delta \le 2.5$, compared to $\delta^2$, which is shown dashed. The $j$ changes at the half-integer points.

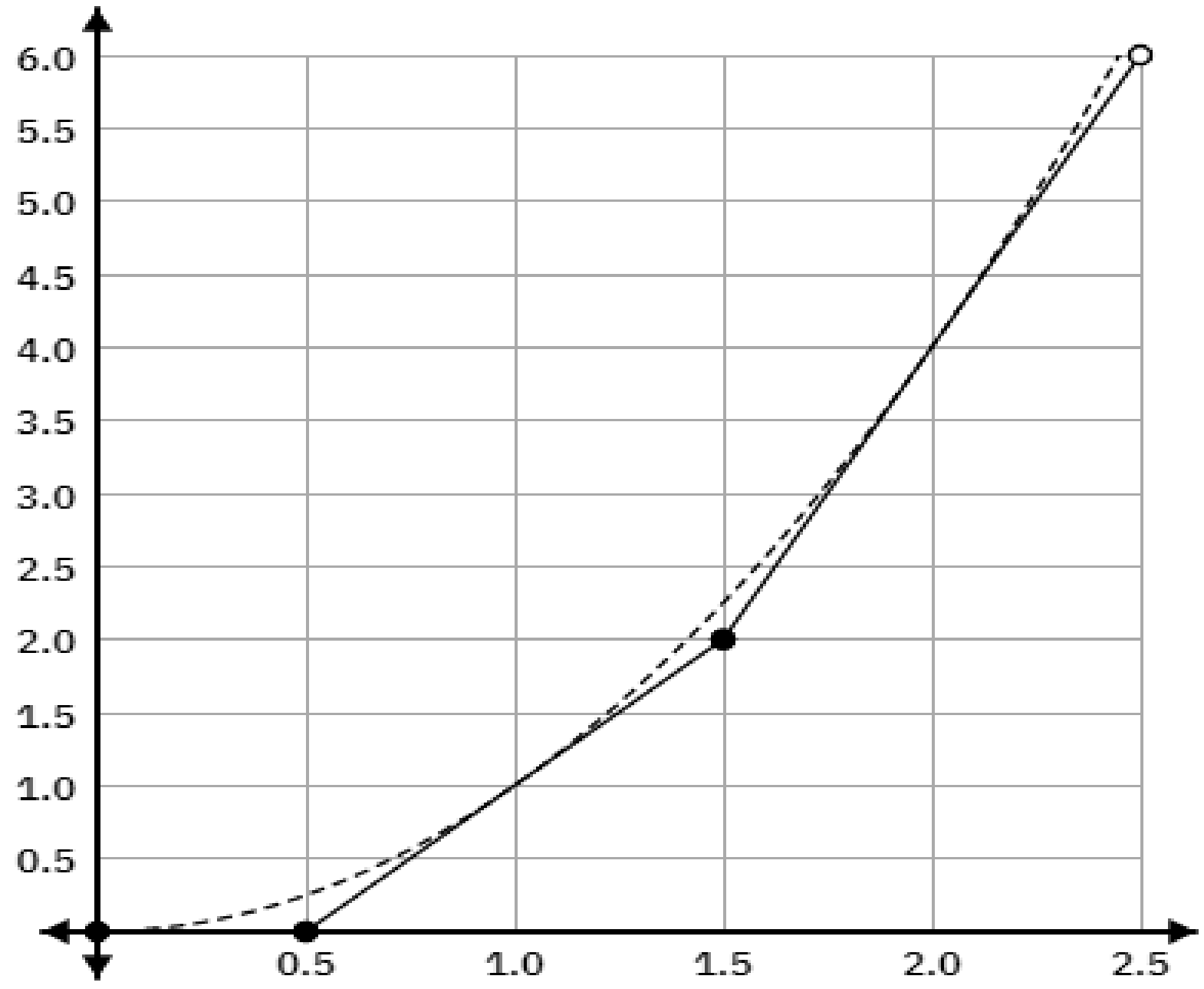

Fig. 7

It is tangent to the parabola at integers, and 1/4 below the parabola at half-integers, with a straight-line graph between half-integer points.

PROOF. By induction. The base case for the induction proof is $j=1$, for which we have already proved it, at (6.4). Assume true for $j$, show true for $j+1$. That is, assume formula is true for $j-1/2 \le \delta \le j+1/2$. Continue to use $n = \left\lfloor \frac{\sqrt{b}}{2} \right\rfloor$. Now consider $j+1/2 \le \delta \le j+3/2$. Break proof into cases:

Case (a) $j+1/2 \le \delta \le j+1$. We need to estimate $V(u+1, b+2m+1)$ for $0 \le m < n$. From (6.2), the distance of $u+1$ from boundary is $d = \frac{\alpha m}{\sqrt{b}} + \delta - 1 + L^{+}\left(\alpha b^{-1/2}/2\right)$, so

$$\frac{\alpha m}{\sqrt{b}} + j - 1/2 \le \frac{\alpha m}{\sqrt{b}} + \delta - 1 + L^{+}\left(\alpha b^{-1/2}/2\right) \le \frac{\alpha m}{\sqrt{b}} + j + L^{+}\left(\alpha b^{-1/2}/2\right).$$

Thus for $0 \le m < n \le \frac{\sqrt{b}}{2}$, $j - 1/2 \le d \le j + \alpha/2 + L^{+}\left(\alpha b^{-1/2}/2\right) < j + 1/2$ since $b \ge b_0$, so the induction hypothesis applies with $d$ in the place of delta, and $b+2m+1$ in the place of $b$. So $V_E(u+1, b+2m+1) = \frac{\alpha}{(b+2m+1)^{3/2}}\left(2j\left(d-1/2\right) - j(j-1) + L(M_j b^{-1/4})\right)$. But $\frac{\alpha}{(b+2m+1)^{3/2}} = \frac{\alpha}{b^{3/2}}\left(1 - L^{+}(1.5b^{-1/2})\right)$, and $2j\left(d-1/2\right) - j(j-1) \le 2j\left(j\right) - j(j-1) = j^2 + j$, so

$$V_E(u+1,b+2m+1)=\frac{\alpha}{b^{3/2}}\begin{pmatrix}2j\left(\frac{\alpha m}{\sqrt{b}}+\delta-3/2+L^+(\alpha b^{-1/2}/2)\right)-j(j-1)\\+L\left(M_jb^{-1/4}+1.5(j^2+j)b^{-1/2}\right)\end{pmatrix}$$

. Calculation with $j\le b^{1/10}$ and $b\ge 1000$ shows $1.5(j^2+j)b^{-1/2}+\alpha jb^{-1/2}\le 2b^{-1/4}$, so

$$V_E(u+1,b+2m+1)=\frac{\alpha}{b^{3/2}}\left(2j\left(\frac{\alpha m}{\sqrt{b}}+\delta-3/2\right)-j(j-1)+L\left((M_j+2)b^{-1/4}\right)\right).$$ So

$$S_{LE}=\frac{\alpha}{b^{3/2}}\sum_{m=0}^{n-1}2^{-2m-1}C_m\left(\left(2j\left(\frac{\alpha m}{\sqrt{b}}+\delta-3/2\right)-j(j-1)+L\left((M_j+2)b^{-1/4}\right)\right)\right)$$

$$=\frac{\alpha}{b^{3/2}}\left(\left(nG_n-1+G_n\right)2j\frac{\alpha}{\sqrt{b}}+\left(2j(\delta-3/2)-j(j-1)\right)\left(1-G_n\right)+L\left((M_j+2)b^{-1/4}\right)\right)$$ . But

$0\le\left(nG_n-1+G_n\right)2j\frac{\alpha}{\sqrt{b}}\le\alpha jb^{-1/4}$, and

$0\ge-\left(2j(\delta-3/2)-j(j-1)\right)G_n\ge-\left(2j(j-1/2)-j(j-1)\right)G_n\ge-j^2b^{-1/4}$, and $j^2\ge\alpha j$, so

$S_{LE}=\frac{\alpha}{b^{3/2}}\left(\left(2j(\delta-3/2)-j(j-1)\right)+L\left(\left(M_j+j^2+2\right)b^{-1/4}\right)\right)$ . From Lemma 6.1, to get $V_E$, add this to $\frac{\alpha}{b^{3/2}}\left(2(\delta-1/2)+L(\delta^2b^{-1/4})\right)$. Noting $\delta\le j+1$, this gives

(6.5) $$V_E(u,b)=\frac{\alpha}{b^{3/2}}\left(\left(2(j+1)(\delta-1/2)-j(j+1)\right)+L\left((M_j+2j^2+2j+3)b^{-1/4}\right)\right).$$

Now we have to deal with

Case(b) $j+1\le\delta\le j+3/2$. Break $m$ into ranges. Let $m_0=\min\left\{\left\lfloor\frac{j+3/2-\delta}{\alpha}\sqrt{b}\right\rfloor,n\right\}$.

Range (i): For $0\le m\le m_0-1,\ \frac{\alpha m}{\sqrt{b}}\le\frac{\alpha m_0}{\sqrt{b}}-\frac{\alpha}{\sqrt{b}}\le j+3/2-\delta-\frac{\alpha}{\sqrt{b}}$ , so

$d=\frac{\alpha m}{\sqrt{b}}+\delta-1+L^+\left(\alpha b^{-1/2}/2\right)\le j+3/2-\delta-\frac{\alpha}{\sqrt{b}}+\delta-1+L^+\left(\alpha b^{-1/2}/2\right)\le j+1/2$, so the induction hypothesis applies with $d$ in the place of delta and $b+2m+1$ in the place of $b$. The same argument that was used in (a) then shows that

$$V_E(u+1,b+2m+1)=\frac{\alpha}{b^{3/2}}\left(2j\left(\frac{\alpha m}{\sqrt{b}}+\delta-3/2\right)-j(j-1)+L\left((M_j+2)b^{-1/4}\right)\right).$$

Range (ii): Suppose $m_0+1 \le m < n$ (this could be empty). Then $m \ge m_0 + 1 \ge \frac{j+3/2-\delta}{\alpha}\sqrt{b}$, so

$d = \frac{\alpha m}{\sqrt{b}} + \delta - 1 + L^+\left(\alpha b^{-1/2}/2\right) \ge j + 3/2 - \delta + \delta - 1 + L^+\left(\alpha b^{-1/2}/2\right) \ge j + 1/2$. And

$\frac{\alpha m}{\sqrt{b}} + \delta - 1 + L^+\left(\alpha b^{-1/2}/2\right) \le \frac{\alpha}{2} + j + \frac{1}{2} + L^+\left(\alpha b^{-1/2}/2\right) \le j+1$ for for $b \ge b_0$. Thus (6.5) applies with $d$ in the place of delta and $b+2m+1$ in the place of $b$, and

$d - 1/2 = \frac{\alpha m}{\sqrt{b}} + \delta - 3/2 + L^+\left(\alpha b^{-1/2}/2\right)$, so

$$V_E(u+1, b+2m+1) = \frac{\alpha}{b^{3/2}}\left(1 - L^+(1.5b^{-1/2})\right)\begin{pmatrix} 2(j+1)\left(\frac{\alpha m}{\sqrt{b}} + \delta - 3/2 + L^+\left(\alpha b^{-1/2}/2\right)\right) - j(j+1) \\ +L\left(\left(M_j + 2j^2 + 2j + 3\right)b^{-1/4}\right)\end{pmatrix}.$$

But $2(j+1)\left(d-1/2\right) - j(j+1) \le 2(j+1)(j+1/2) - j(j+1) = j^2+2j+1$, and computation shows $1.5(j^2+2j+1)b^{-1/2} + \alpha(j+1)b^{-1/2} \le 3b^{-1/4}$ for $b \ge b_0$, so

$$V_E(u+1, b+2m+1) = \frac{\alpha}{b^{3/2}}\left(2(j+1)\left(\frac{\alpha m}{\sqrt{b}} + \delta - 3/2\right) - j(j+1) + L\left(\left(M_j + 2j^2 + 2j + 6\right)b^{-1/4}\right)\right).$$

Range (iii): For $m = m_0$, either $d \le 1/2$ or $d \ge 1/2$, so one of the two formulas applies. But

$m_0 = \frac{j+3/2-\delta}{\alpha}\sqrt{b} - f, 0 \le f < 1$, so $\frac{\alpha m_0}{\sqrt{b}} - 3/2 + \delta = j - \frac{\alpha f}{\sqrt{b}}$, so

$2j\left(\frac{\alpha m_0}{\sqrt{b}} + \delta - 3/2\right) - j(j-1) = j^2 + j - 2j\alpha f b^{-1/2}$ and

$2(j+1)\left(\frac{\alpha m_0}{\sqrt{b}} + \delta - 3/2\right) - j(j+1) = j^2 + j - 2(j+1)\alpha f b^{-1/2}$, which differ by only $2\alpha f b^{-1/2}$. We can use the formula from (ii) for $m = m_0$ regardless, and the result is the same, since

$(M_j+2)b^{-1/4} + 2\alpha f b^{-1/2} = L\left(\left(M_j + 2j^2 + 2j + 6\right)b^{-1/4}\right)$.

Summing over the formulas for ranges (i)-(iii),

$$S_{LE} = \frac{\alpha}{b^{3/2}}\left\{\begin{matrix} \sum_{m=0}^{m_0-1} 2^{-2m-1}C_m\left(2j\left(\frac{\alpha m}{\sqrt{b}} + \delta - 3/2\right) - j(j-1) + L\left((M_j+2)b^{-1/4}\right)\right) \\ + \sum_{m=m_0}^{n-1} 2^{-2m-1}C_m\begin{pmatrix} 2(j+1)\left(\frac{\alpha m}{\sqrt{b}} + \delta - 3/2\right) - j(j+1) \\ +L\left(\left(M_j + 2j^2 + 2j + 6\right)b^{-1/4}\right)\end{pmatrix}\end{matrix}\right\}.$$

Now $\sum_{m=0}^{n} 2^{-2m-1} C_m < 1$, and $(M_j + 2)b^{-1/4} = L\left(\left(M_j + 2j^2 + 2j + 6\right)b^{-1/4}\right)$, so we can get the error terms out of the summation.

$$S_{LE} = \frac{\alpha}{b^{3/2}} \begin{Bmatrix} \sum_{m=0}^{m_0-1} 2^{-2m-1} C_m \left( 2j\left(\frac{\alpha m}{\sqrt{b}} + \delta - 3/2\right) - j(j-1) \right) \\ + \sum_{m=m_0}^{n-1} 2^{-2m-1} C_m \left( 2(j+1)\left(\frac{\alpha m}{\sqrt{b}} + \delta - 3/2\right) - j(j+1) \right) \\ + L\left(\left(M_j + 2j^2 + 2j + 6\right)b^{-1/4}\right) \end{Bmatrix}$$

$$= \frac{\alpha}{b^{3/2}} \begin{Bmatrix} \sum_{m=0}^{n-1} 2^{-2m-1} C_m \left( 2j\left(\frac{\alpha m}{\sqrt{b}} + \delta - 3/2\right) - j(j-1) \right) \\ + \sum_{m=m_0}^{n-1} 2^{-2m-1} C_m \left( 2\left(\frac{\alpha m}{\sqrt{b}} + \delta - 3/2\right) - 2j \right) + L\left(\left(M_j + 2j^2 + 2j + 6\right)b^{-1/4}\right) \end{Bmatrix}$$

$$= \frac{\alpha}{b^{3/2}} \begin{Bmatrix} 2j\frac{\alpha}{\sqrt{b}}\left(nG_n - 1 + G_n\right) + \left(2j(\delta - 1/2) - j(j+1)\right)\left(1 - G_n\right) + \frac{2\alpha}{\sqrt{b}} \sum_{m=m_0}^{n-1} 2^{-2m-1} C_m m \\ + 2\left(\delta - j - 3/2\right)(G_{m_0} - G_n) + L\left(\left(M_j + 2j^2 + 2j + 6\right)b^{-1/4}\right) \end{Bmatrix}.$$

Now to estimate some terms.

$2j\frac{\alpha}{\sqrt{b}}\left(nG_n - 1 + G_n\right) \le 2j\frac{\alpha}{n/2} nG_n \le \alpha j\sqrt{2.14/\pi}\, b^{-1/4}$.

$\left(2j(\delta - 1/2) - j(j+1)\right)G_n \le (j^2 + j)G_n \le (j^2 + j)\sqrt{2.14/\pi}\, b^{-1/4}$.

$\frac{2\alpha}{\sqrt{b}} \sum_{m=m_0}^{n-1} 2^{-2m-1} C_m m \le \frac{2\alpha}{\sqrt{b}} \sum_{m=0}^{n-1} 2^{-2m-1} C_m m \le \frac{2\alpha}{\sqrt{b}}(nG_n) \le \alpha\sqrt{2/\pi}\, b^{-1/4}$.

$2\left(j + 3/2 - \delta\right)(G_{m_0} - G_n) \le 2\left(j + 3/2 - \delta\right)G_{m_0} \le 2\frac{\alpha}{\sqrt{b}}\left(m_0 + 1\right)G_{m_0} \le 2\frac{\alpha}{\sqrt{b}} nG_n + 2\frac{\alpha}{\sqrt{b}}$, because

$mG_m$ is an increasing function and $G_n$ is decreasing. Thus

$2\left(j + 3/2 - \delta\right)(G_{m_0} - G_n) \le \alpha\sqrt{2/\pi}\, b^{-1/4} + 2\alpha b^{-1/2} = L\left(b^{-1/4}\right)$.

Putting it all together, $S_{LE} = \frac{\alpha}{b^{3/2}}\left\{2j(\delta - 1/2) - j(j+1) + L\left(\left(M_j + 3j^2 + 4j + 8\right)b^{-1/4}\right)\right\}$ is good enough. Add this to $\frac{\alpha}{b^{3/2}}\left(2(\delta - 1/2) + L(\delta^2 b^{-1/4})\right)$ from Lemma 6.1, using $\delta \le j + 3/2$, to get

$V_E(u,b) = \frac{\alpha}{b^{3/2}}\left(\left(2(j+1)(\delta - 1/2) - j(j+1)\right) + L\left(\left(M_j + 4j^2 + 7j + 11\right)b^{-1/4}\right)\right)$,

valid over the range $j + 1/2 \le \delta \le j + 3/2$, and the induction proof is complete, upon letting

$M_{j+1} = M_j + 4j^2 + 7j + 11$. This recursion, with $M_1 = 5$, is easily seen to imply $M_j \le 5j^3$. Of course a much smaller upper bound on $M_j$ is possible if we cared, but we don't. □

**7. Proof of Theorem 1.3.** In this section, let $n_0 = \left\lfloor \frac{\sqrt{b}}{\alpha} \right\rfloor$, so $\frac{\alpha n_0}{\sqrt{b}} = 1 - L^+(\alpha b^{-1/2})$ and $\frac{\alpha}{\sqrt{b}} = \frac{1}{n_0} - L^+(\alpha b^{-1})$ ; we'll see that replacing $\frac{\sqrt{b}}{\alpha}$ by its integer part isn't going to matter, to the order of interest. Let $0 < p \le 1/10$, to be decided later (in the end, it will be 1/12). Let $J + 1 = \lfloor b^p \rfloor$, so $J \sim b^p$ . In this section, we are going to let $n = (J+1)n_0$, so $n$ will go to infinity faster than $\sqrt{b}$ , but only slightly. This will larger $n$ will cause some of the leaf Values to be more than just the ratio. By dividing $n$ into $J$ stretches of size $n_0$, we'll get a linear approximation to the extra part of the leaf Values, $V_E(u+1, b+2m+1)$, on each stretch, via Theorem 2; that's Lemma 7.1. As usual, $\delta = \alpha\sqrt{b} - u$ . Assume throughout this section that $b \ge b_0$. We may then assume $\frac{1}{2} - \frac{.44}{\sqrt{\pi}} b^{-1/4} \le \delta \le \frac{1}{2} - \frac{.13}{\sqrt{\pi}} b^{-1/4}$ because we already know that the exact stop value occurs for some delta in this range. Say $\delta = \frac{1}{2} - \gamma b^{-1/4}$, where $\frac{.13}{\sqrt{\pi}} \le \gamma \le \frac{.44}{\sqrt{\pi}}$. By letting $J$ go to infinity with $b$, slowly, at just the right rate, we can make the upper and lower bounds from Lemma 5.2 come together, to order $o(b^{-1/4})$, as we'll later see.

LEMMA 7.1. There exists $K$ such that for $\delta = 1/2 - L^+(b^{-1/4})$ and for $j = 1, \ldots, J$,

$$V_E(u+1, b+2m+1) = \frac{\alpha}{b^{3/2}}\left(2j\frac{m}{n_0} - j(j+1) + L\left(Kj^3 b^{-1/4}\right)\right),\ jn_0 \le m \le (j+1)n_0 .$$

PROOF: Let $jn_0 \le m \le (j+1)n_0$. Then

$$d = \alpha\sqrt{b+2m+1} - (u+1) = \alpha\sqrt{b}\left(b+2m+1\right)^{1/2} - (\alpha\sqrt{b} - \delta - 1) .$$

Now

$$\left(1 + \frac{2m+1}{b}\right)^{1/2} \le 1 + \frac{1}{2}\frac{2m+1}{b} - \frac{1}{8}\frac{(2m+1)^2}{b^2} + \frac{1}{16}\frac{(2m+1)^3}{b^3} = 1 + \frac{m}{b} + \frac{1}{2b} - \frac{1}{8}\frac{(2m+1)^2}{b^2}\left(1 - \frac{2m+1}{2b}\right)$$

$$\le 1 + \frac{m}{b} + \frac{1}{2b} - \frac{1}{2}\frac{{n_0}^2}{b^2}\left(1 - \frac{n_0}{b}\right)$$

since $x^2\left(1-x/b\right)$ is decreasing over the range of $x$ here. Now $n_0 \geq \frac{\sqrt{b}}{\alpha}-1$, and putting that in, $\frac{1}{2}\frac{n_0^{\;2}}{b^2}\left(1-\frac{n_0}{b}\right)\cong\frac{1}{2\alpha^2 b}$ , a little smaller, but close enough when $b \geq b_0$ so that $\frac{1}{2b}-\frac{1}{2}\frac{n_0^{\;2}}{b^2}\left(1-\frac{n_0}{b}\right)<0$ , as can be shown. Thus $\left(1+\frac{2m+1}{b}\right)^{1/2}\leq 1+\frac{m}{b}$, so

$d=\alpha\sqrt{b}\left(1+\frac{2m+1}{b}\right)^{1/2}-\alpha\sqrt{b}+\delta-1\leq\frac{\alpha m}{\sqrt{b}}+\delta-1\leq\frac{m}{n_0}-\frac{1}{2}-\gamma b^{-1/4}$ . Going the other way,

$$\left(1+\frac{2m+1}{b}\right)^{1/2}\geq 1+\frac{m}{b}+\frac{1}{2b}-\frac{1}{8}\frac{(2m+1)^2}{b^2}\geq 1+\frac{m}{b}+\frac{1}{2b}-\frac{1}{8}\frac{\left(2b^p\sqrt{b}/\alpha+1\right)^2}{b^2}\geq 1+\frac{m}{b}-\frac{1}{2\alpha^2}b^{-1+2p}.$$

So $d\geq\frac{\alpha m}{\sqrt{b}}-\frac{1}{2\alpha}b^{-1/2+2p}+\delta-1\geq\frac{m}{n_0}-m\alpha b^{-1}-\frac{1}{2\alpha}b^{-1/2+2p}-\frac{1}{2}-\gamma b^{-1/4}$

$$\geq\frac{m}{n_0}-\frac{1}{2}-b^{-1/2+p}-\frac{1}{2\alpha}b^{-1/2+2p}-\gamma b^{-1/4}\geq\frac{m}{n_0}-\frac{1}{2}-b^{-1/4}.$$

Putting together, $d=\frac{m}{n_0}-\frac{1}{2}-\psi b^{-1/4}$, for some $0\leq\psi\leq 1$, is good enough.

Note $j-1/2\leq\frac{m}{n_0}-1/2\leq j+1/2$.

(i) If $j-1/2\leq d=\frac{m}{n_0}-1/2-\psi b^{-1/4}\leq j+1/2$, then Theorem 6.2 gives

$$V_E(u+1,b+2m+1)=\frac{\alpha}{(b+2m+1)^{3/2}}\left(2j\left(\frac{m}{n_0}-1/2-\psi b^{-1/4}-1/2\right)-j(j-1)+L\left(5j^3b^{-1/4}\right)\right).$$

Now $\frac{\alpha}{(b+2m+1)^{3/2}}=\frac{\alpha}{b^{3/2}}\left(1-L^+\left(1.5(2m+1)b^{-1}\right)\right)$. But $m\leq(J+1)n_0\leq b^{1/10}b^{-1/2}/\alpha$, so calculation shows $1.5(2m+1)b^{-1}\leq 1.27b^{-1/4}$ for $b\geq b_0\geq 1000$. Then

$$V_E(u+1,b+2m+1)=\frac{\alpha}{b^{3/2}}\left(1-L^+(1.27b^{-1/4})\right)\left(2j\frac{m}{n_0}-j(j+1)+L\left((5j^3+2j)b^{-1/4}\right)\right)$$

$=\frac{\alpha}{b^{3/2}}\left(2j\frac{m}{n_0}-j(j+1)+L\left((5j^3+2j^2+2j)b^{-1/4}\right)\right)$, since $2j\frac{m}{n_0}-j(j+1)\leq j^2$.

(ii) Suppose, however, that $d=\frac{m}{n_0}-1/2-\psi b^{-1/4}\leq j-1/2$: the $\psi b^{-1/4}$ term bumped us down into the next interval below. But just barely. $\frac{m}{n_0}-1/2\geq j-1/2$, and $\frac{m}{n_0}-1/2-\psi b^{-1/4}\leq j-1/2$, so

$\frac{m}{n_0} = j + L^+(\psi b^{-1/4})$. We compute using Theorem 6.2, but with $j-1$ in place of $j$, so

$V_E(u+1, b+2m+1)$

$$= \frac{\alpha}{(b+2m+1)^{3/2}}\left(2(j-1)\left(\frac{m}{n_0} - 1/2 - \psi b^{-1/4} - 1/2\right) - (j-1)(j-2) + L\left(5(j-1)^3 b^{-1/4}\right)\right)$$

$$= \frac{\alpha}{(b+2m+1)^{3/2}}\left(2j\frac{m}{n_0} - 2\frac{m}{n_0} - 2(j-1)\left(1+\psi b^{-1/4}\right) - (j-1)(j-2) + L\left(5(j-1)^3 b^{-1/4}\right)\right)$$

$$= \frac{\alpha}{(b+2m+1)^{3/2}}\left(2j\frac{m}{n_0} - 2j - 2L^+(\psi b^{-1/4}) - 2(j-1) - 2(j-1)\psi b^{-1/4} - (j-1)(j-2) + L\left(5(j-1)^3 b^{-1/4}\right)\right)$$

$$= \frac{\alpha}{(b+2m+1)^{3/2}}\left(2j\frac{m}{n_0} - j(j+1) - 2L^+(\psi b^{-1/4}) - 2j\psi b^{-1/4} + 2\psi b^{-1/4} + L\left(5(j-1)^3 b^{-1/4}\right)\right)$$

$$= \frac{\alpha}{(b+2m+1)^{3/2}}\left(2j\frac{m}{n_0} - j(j+1) + L\left(\left(5(j-1)^3 + 2j\right) b^{-1/4}\right)\right).$$ Then the same argument as in (i) gives

$$V_E(u+1, b+2m+1) = \frac{\alpha}{b^{3/2}}\left(2j\frac{m}{n_0} - j(j+1) + L\left(\left(5(j-1)^3 + 2j^2 + 2j\right) b^{-1/4}\right)\right)$$

Cases (i) and (ii) are both covered by $V_E(u+1, b+2m+1)$

$= \frac{\alpha}{b^{3/2}}\left(2j\frac{m}{n_0} - j(j+1) + L\left(9j^3 b^{-1/4}\right)\right)$, so $K = 9$ is good enough. □

To prove the main theorem, we go down to row $(J+1)n_0$ in the backwards induction tree. We have to put in the extra contribution $S_{LE}$ to the leaf sum that comes from the leaf values that exceed the ratio.

LEMMA 7.2. $S_{LE} = \frac{\alpha}{b^{3/2}}\left\{C_{LE} G_n + O\left(J^{5/2} b^{-1/2}\right)\right\}$, where $C_{LE} = \frac{1}{3}J^2 - \frac{1}{3}J - 4\zeta(-1/2)\sqrt{J+1} - \frac{5}{6} + O(J^{-1})$. $J$ and $n$ are as asserted at the top of this section.

PROOF: Let $S_{LE}(j)$ denote the contribution to $S_{LE}$ from (7.1) for $m$ in the range $jn_0 \le m < (j+1)n_0$. $S_{LE}(j)$ is 0 for $j = 0$. Assume now $j \ge 1$. From Lemma 7.1,

$$S_{LE}(j) = \frac{\alpha}{b^{3/2}}\left\{\frac{2j}{n_0}\sum_{m=jn_0}^{(j+1)n_0-1} m2^{-2m-1}C_m - \left(j(j+1) + L\left(Kj^3 b^{-1/4}\right)\right)\sum_{m=jn_0}^{(j+1)n_0-1} 2^{-2m-1}C_m\right\}$$

$$=\frac{\alpha}{b^{3/2}}\left\{\frac{2j}{n_0}\left((j+1)n_0G_{(j+1)n_0}-jn_0G_{jn_0}+G_{(j+1)n_0}-G_{jn_0}\right)-\left(j(j+1)+L\left(Kj^3b^{-1/4}\right)\left(G_{jn_0}-G_{(j+1)n_0}\right)\right)\right\}$$

$$=\frac{\alpha}{b^{3/2}}\left\{3j(j+1)G_{(j+1)n_0}-j(3j+1)G_{jn_0}+\frac{2j}{n_0}\left(G_{(j+1)n_0}-G_{jn_0}\right)-\left(L\left(Kj^3b^{-1/4}\right)\left(G_{jn_0}-G_{(j+1)n_0}\right)\right)\right\}.$$

But

$$G_{jn_0}-G_{(j+1)n_0}\le\frac{1}{\sqrt{\pi}}\left(\frac{1}{\sqrt{n_0j}}-\frac{1}{\sqrt{n_0(j+1)+1/2}}\right)=\frac{1}{\sqrt{\pi}}\left(\frac{\sqrt{n_0(j+1)+1/2}-\sqrt{n_0j}}{\sqrt{n_0(j+1)+1/2}\sqrt{n_0j}}\right)$$

$$=\frac{1}{\sqrt{\pi}}\left(\frac{n_0+1/2}{\sqrt{n_0(j+1)+1/2}\sqrt{n_0j}\left(\sqrt{n_0(j+1)+1/2}+\sqrt{n_0j}\right)}\right)\le\frac{1}{\sqrt{\pi n_0}}\left(\frac{1}{2j^{3/2}}\right),$$ so

$$S_{LE}(j)=\frac{\alpha}{b^{3/2}}\left\{3j(j+1)G_{(j+1)n_0}-j(3j+1)G_{jn_0}+L\left(Kj^{3/2}b^{-1/2}\right)\right\}.$$

Note in general $\frac{1}{\sqrt{n\pi}}\left(1-\frac{1}{4n}\right)\le G_n\le\frac{1}{\sqrt{n\pi}}$. So $\frac{G_{jn_0}}{G_{(J+1)n_0}}=\frac{\sqrt{J+1}}{\sqrt{j}}\left(1-L\left(\frac{1}{4jn_0}\right)\right)$, and

$$S_{LE}(j)=\frac{\alpha}{b^{3/2}}\left\{\begin{array}{l}\left[3j(j+1)\frac{\sqrt{J+1}}{\sqrt{j+1}}\left(1+L\left(\frac{1}{4(j+1)n_0}\right)\right)-j(3j+1)\frac{\sqrt{J+1}}{\sqrt{j}}\left(1+L\left(\frac{1}{4jn_0}\right)\right)\right]G_{(J+1)n_0}\\+L\left(Kj^{3/2}b^{-1/2}\right)\end{array}\right\}$$

$$=\frac{\alpha}{b^{3/2}}\left\{\left[3j\sqrt{(j+1)}-\sqrt{j}(3j+1)\right]\sqrt{J+1}G_{(J+1)n_0}+L\left((K+1)j^{3/2}b^{-1/2}\right)\right\}$$

$$=\frac{\alpha}{b^{3/2}}\left\{\left[3j\sqrt{(j+1)}-3(j-1)\sqrt{j}-4\sqrt{j}\right]\sqrt{J+1}G_{(J+1)n_0}+L\left((K+1)j^{3/2}b^{-1/2}\right)\right\}.$$

We wrote it that way so that telescoping occurs for the first two terms when summing:

$$S_{LE}=\frac{\alpha}{b^{3/2}}\left\{\left[\sum_{j=1}^{J}3j\sqrt{(j+1)}-3(j-1)\sqrt{j}-4\sum_{j=1}^{J}\sqrt{j}\right]\sqrt{J+1}G_{(J+1)n_0}+L\left((K+1)J^{5/2}b^{-1/2}\right)\right\}$$

$$=\frac{\alpha}{b^{3/2}}\left\{\left[3J\sqrt{(J+1)}-4\sum_{j=1}^{J}\sqrt{j}\right]\sqrt{J+1}G_n+O\left(J^{5/2}b^{-1/2}\right)\right\}.$$ At this point we are through summing over an unbounded range, so we replaced the big-L with big-O notation, with no harm.

There is a well-known asymptotic formula for the sum of square roots as a generalized harmonic number: $H_J^{(-1/2)}=\sum_{j=1}^{J}\sqrt{j}=\frac{2J^{3/2}}{3}+\frac{J^{1/2}}{2}+\zeta(-1/2)+\frac{J^{-1/2}}{24}+O\left(J^{-5/2}\right)$ [7, pg. 594, prob. 9.27].

$$S_{LE}=\frac{\alpha}{b^{3/2}}\left\{\left[3J\sqrt{(J+1)}-4\left(\frac{2J^{3/2}}{3}+\frac{J^{1/2}}{2}+\zeta(-1/2)+\frac{J^{-1/2}}{24}+O\left(J^{-5/2}\right)\right)\right]\sqrt{J+1}G_{(J+1)n_0}+O\left(J^{5/2}b^{-1/2}\right)\right\}$$

$$=\frac{\alpha}{b^{3/2}}\left\{\begin{array}{l}\left[3J(J+1)-4\left(\frac{2J^2}{3}+\frac{J}{2}+\frac{1}{24}+O\left(J^{-5/2}\right)\right)\left(1+\frac{1}{J}\right)^{1/2}-4\zeta(-1/2)\sqrt{J+1}\right]G_{(J+1)n_0}\\+O\left(J^{5/2}b^{-1/2}\right)\end{array}\right\}$$

$$=\frac{\alpha}{b^{3/2}}\left\{\begin{array}{l}\left[3J(J+1)-4\left(\frac{2J^2}{3}+\frac{J}{2}+\frac{1}{24}+O\left(J^{-5/2}\right)\right)\left(1+\frac{1}{2J}-\frac{1}{8J^2}+O\left(J^{-3}\right)\right)-4\zeta(-1/2)\sqrt{J+1}\right]G_n\\+O\left(J^{5/2}b^{-1/2}\right)\end{array}\right\}$$

$$=\frac{\alpha}{b^{3/2}}\left\{\left[\frac{1}{3}J^2-\frac{1}{3}J-\frac{5}{6}+O(J^{-1})-4\zeta(-1/2)\sqrt{J+1}\right]G_n+O\left(J^{5/2}b^{-1/2}\right)\right\}$$ . Let

$$C_{LE}=\frac{1}{3}J^2-\frac{1}{3}J-\frac{5}{6}+O(J^{-1})-4\zeta(-1/2)\sqrt{J+1}\,.$$ □

Finally, to prove Theorem 1.3, go back to Lemma 5.2, to get the bounds on the Tree Sum using this $n$. Our assumption about $J$ implies $n=O(b^{1/2+p})$ with $p\le 1/10$, so Lemma 5.2 applies. For the upper bound, $TreeSum(n,u,b)-\frac{u}{b}\le\frac{\alpha}{b^{3/2}}\left\{C+\delta B+\delta^2 A\right\}+S_{LE}=\frac{\alpha}{b^{3/2}}\left\{C^*+\delta B+\delta^2 A\right\}$, where $C^*=-\left(1-(C_1+C_{LE})G_n\right)+O(J^{5/2}b^{-1/2})$, and $C_1=1+2\frac{\alpha n}{\sqrt{b}}-\frac{\alpha^2}{3}\frac{n^2}{b}$, $A=G_n$, $B=2\left(1-B_1G_n\right)$ with $B_1=1+\alpha\frac{n}{\sqrt{b}}$. Solving the quadratic equation just as we did in proving Lemma 5.3,

$$\delta_0=\frac{-B+B\left(1-\frac{4AC^*}{B^2}\right)^{1/2}}{2A}=\frac{-B+B\left(1-\frac{2AC^*}{B^2}-\frac{2A^2C^{*2}}{B^4}-O(A^3)\right)}{2A}=-\frac{C^*}{B}-\frac{AC^{*2}}{B^3}+O(b^{-1/2-p}),$$

since $A=G_n=O(b^{-1/4-p/2})$. Now $-\frac{C^*}{B}-\frac{AC^{*2}}{B^3}$

$$=\frac{1-(C_1+C_{LE})G_n+O(J^{5/2}b^{-1/2})}{2(1-B_1G_n)}-G_n\frac{\left(1-(C_1+C_{LE})G_n+O(J^{5/2}b^{-1/2})\right)^2}{8\left(1-B_1G_n\right)^3}$$

$$=\frac{1}{2}\left(1-(C_1+C_{LE})G_n\right)\left(1+B_1G_n+O(J^{5/2}b^{-1/2})\right)-\frac{G_n}{8}+O(J^{5/2}b^{-1/2})$$

$$=\frac{1}{2}-\frac{1}{2}\left(C_1+C_{LE}-B_1+\frac{1}{4}\right)G_n+O(J^{5/2}b^{-1/2})=\frac{1}{2}-\frac{1}{2}\left(C_{LE}+\frac{\alpha n}{\sqrt{b}}-\frac{1}{3}\frac{\alpha^2n^2}{b}+\frac{1}{4}\right)G_n+O(J^{5/2}b^{-1/2})\,.$$

But $\frac{\alpha n}{\sqrt{b}} = J+1+O\left((J+1)b^{-1/2}\right)$, so $\frac{\alpha n}{\sqrt{b}} - \frac{1}{3}\frac{\alpha^2 n^2}{b} + \frac{1}{4} = J+1-\frac{1}{3}(J+1)^2 + \frac{1}{4} + O\left((J+1)^2 b^{-1/2}\right)$

$= -\frac{1}{3}J^2 + \frac{1}{3}J + \frac{11}{12} + O\left((J+1)^2 b^{-1/2}\right)$, so $C_{LE} + \frac{\alpha n}{\sqrt{b}} - \frac{1}{3}\frac{\alpha^2 n^2}{b} + \frac{1}{4}$

$= \frac{1}{3}J^2 - \frac{1}{3}J - \frac{5}{6} + O(J^{-1}) - 4\zeta(-1/2)\sqrt{J+1} - \frac{1}{3}J^2 + \frac{1}{3}J + \frac{11}{12} + O\left((J+1)^2 b^{-1/2}\right)$

$= -4\zeta(-1/2)\sqrt{J+1} + \frac{1}{12} + O(J^{-1})$. The higher power terms miraculously canceled!, exposing the zeta function term as dominant.

$$\delta_0 = \frac{1}{2} - \frac{1}{2}\left(-4\zeta(-1/2)\sqrt{J+1} + \frac{1}{12} + O(J^{-1})\right)\frac{1}{\sqrt{\pi(J+1)n_0}} + O\left(J^{5/2}b^{-1/2}\right)$$

$$= \frac{1}{2} - \frac{1}{2}\left(-4\zeta(-1/2) + \frac{1}{12\sqrt{J+1}} + O(J^{-1})\right)\frac{\sqrt{\alpha}}{\sqrt{\pi}}b^{-1/4} + O(J^{5/2}b^{-1/2}) \qquad (*)$$

$$= \frac{1}{2} - \frac{\left(-2\zeta(-1/2)\right)\sqrt{\alpha}}{\sqrt{\pi}}b^{-1/4} + O(J^{-1/2}b^{-1/4}) + O(J^{5/2}b^{-1/2}).$$

Since $J \sim b^p$, the best we can do with this is to let $p = 1/12$. We get

$$\delta_0 = \frac{1}{2} - \frac{\left(-2\zeta(-1/2)\right)\sqrt{\alpha}}{\sqrt{\pi}}b^{-1/4} + O(b^{-7/24}) = \frac{1}{2} - \frac{\left(-2\zeta(-1/2)\right)\sqrt{\alpha}}{\sqrt{\pi}}b^{-1/4} + o(b^{-1/4}).$$

For the lower bound, to get $\delta_0'$, the only change is that $C_1 = 1 + 2\frac{\alpha n}{\sqrt{b}} - \frac{\alpha^2}{3}\frac{n^2}{b}$ is replaced by $C_1' = 7/12 + 2\frac{\alpha n}{\sqrt{b}} - \frac{\alpha^2}{3}\frac{n^2}{b}$, the effect of which is to change the $\frac{1}{12\sqrt{J+1}}$ term in (*) to $-\frac{1}{3\sqrt{J+1}}$, so that $\delta_0'$ is the same asymptotically as $\delta_0$, to order $O(b^{-7/24})$. This completes the proof of Theorem 1.1. □

As a final remark, we could say that we expected those higher terms to cancel, based on the idea that going further down the tree leads to less weight on the row and more weight on the leaves, where, at least for a while, the errors are quite small thanks to Theorem 6.2. But honestly, when that happened with just the right choice of $p$, we thanked Tyche rather than crediting our insight, since we didn't really know if it would happen before doing it.

**Acknowledgements**. We thank T. P. Hill for introducing us to this problem, and for pointing out the work of Häggström and Wästlund, which motivated what we did. We also thank him for discussions and encouragement about what we were doing, which was much needed lest we give up when the algebra seemed hopeless. We thank Doron Zeilberger for encouragement and inspiration. We thank Steve Mann for supplying the fast Intel MMX code for both floating and integer 128-bit arithmetic that was needed for the rigorous determination of stop rules for the extremely large number of flips considered.